\documentclass{amsart}
\usepackage{amsmath,amssymb,amsthm,mathrsfs}

\newcommand{\sF}{{\mathcal F}}

\newcommand{\sW}{{\mathcal W}}

\def\a{{\alpha}}
\def\b{\beta}
\def\c{\gamma}
\def\d{\delta}
\def\de{\Delta}

\def\l{\lambda}

\def\om{\Omega}

\def\va{\varphi}
\def\vp{\varepsilon}

\def\z{\zeta}
\def\ts{\times}

\def\iy{\infty}

\def\BC{{\mathbb C}}

\def\BT{{\mathbb T}}

\renewcommand{\theequation}{\arabic{section}.\arabic{equation}}

\newcommand{\bpr}{{\noindent\textbf{Proof.}\ \ }}
\newcommand{\epr}{{$\mbox{}$ \hfill $\Box$}}

\newcommand{\mat}[2]{\ensuremath{\left[\begin{array}{#1}#2\end{array}\right]}}

\newtheorem{thm}{Theorem}[section]
\newtheorem{prop}[thm]{Proposition}
\newtheorem{lem}[thm]{Lemma}

\newtheorem{rem}[thm]{Remark}

\newcommand{\ands}{\quad\mbox{and}\quad}

\begin{document}

\title{The discrete twofold Ellis-Gohberg inverse problem}

\author[S. ter Horst]{S. Ter Horst}
\address{S. Ter Horst, Department of Mathematics, Unit for BMI, North-West
University,
Potchefstroom, 2531 South Africa}
\email{Sanne.TerHorst@nwu.ac.za}

\author[M.A. Kaashoek]{M.A. Kaashoek}
\address{M.A. Kaashoek, Department of Mathematics,
VU University Amsterdam\\
De Boelelaan 1081a, 1081 HV Amsterdam, The Netherlands}
\email{m.a.kaashoek@vu.nl}

\author[F. van Schagen]{F. van Schagen}
\address{F. van Schagen, Department of Mathematics,
VU University Amsterdam\\
De Boelelaan 1081a, 1081 HV Amsterdam, The Netherlands}
\email{f.van.schagen@vu.nl}

\thanks{This work is based on the research supported in part by the
National Research Foundation of South Africa (Grant Number 93406).}

\begin{abstract}
In this paper a twofold inverse problem for orthogonal matrix functions in the Wiener class is considered. The scalar-valued version of this problem was solved by Ellis and Gohberg in 1992. Under reasonable conditions, the problem is reduced to an invertibility condition on an operator that is defined using the Hankel and Toeplitz operators associated to the Wiener class functions that comprise the data set of the inverse problem. It is also shown that in this case the solution is unique. Special attention is given to the case that the Hankel operator of the solution is a strict contraction and the case where the functions are matrix polynomials.
\end{abstract}

\subjclass[2010]{Primary: 47A56; Secondary: 47B35, 47A50, 15A29}

\keywords{Inverse problem, Wiener algebra, Toeplitz operator, Hankel operator, structured operators, operator inversion}

\maketitle


\section{Introduction}
To state  our  main problem we need  some notation and  terminology about Wiener class functions.
Throughout  $\sW^{n\ts m} $ denotes the space of $ n\ts m$ matrix functions with entries in the Wiener algebra on the unit circle. Thus a matrix function  $\va $ belongs to $\sW^{n\ts m} $  if and only if $\va$  is continuous on the unit circle and   its Fourier coefficients $\dots \va_{-1}, \va_0, \va_1, \ldots$ are absolutely summable. We set
\begin{align*}
\sW_+^{n\ts m}&=\{\va\in \sW^{n\ts m}\mid \va_j= 0, \quad
\mbox{for } j=-1, -2, \dots\},   \\
\sW_-^{n\ts m}&=\{\va\in \sW^{n\ts m}\mid \va_j= 0, \quad \mbox{for } j=1, 2, \dots\}, \\
\sW_d^{n\ts m} &=\{\va\in \sW^{n\ts m}\mid \va_j= 0, \quad \mbox{for } j\not =0\}, \\
\sW_{+,0}^{n\ts m}&=\{\va\in \sW^{n\ts m}\mid \va_j= 0, \quad \mbox{for } j=0,-1, -2, \dots\},   \\
\sW_{-,0}^{n\ts m}&=\{\va\in \sW^{n\ts m}\mid \va_j= 0, \quad \mbox{for } j=0,1, 2, \dots\}.
\end{align*}
Given $\va \in \sW^{n\ts m}$ the function $\va^*$ is defined by $\va^*(\z)=\va(\z)^*$ for each $\z\in \BT$.  Thus the $j$-th Fourier coefficient of  $\va^*$ is given by $(\va^*)_j=(\va_{-j})^*$.  The map $\va \mapsto  \va^*$ defines an involution  which transforms  $\sW^{n\ts m}$ into $\sW^{m \times n}$,   $ \sW_+^{n\ts m}$ into $ \sW_-^{m \ts n}$, $ \sW_{-, 0}^{n\ts m}$ into $\sW_{+,0}^{m\ts n}$, etc.

The data of the inverse problem we shall be dealing with consist of four functions, namely
\begin{equation}\label{Fabcd2}
\a\in \sW_+^{p\ts p}, \quad \b\in \sW_+^{p\ts q}, \quad \c\in \sW_-^{q\ts p}, \quad \d\in \sW_-^{q\ts q},
\end{equation}
and we are interested in finding  $g\in \sW_+^{p\ts q}$ such that
\begin{align}
& \a + g \c - e_p \in  \sW_{-,0}^{ p \times p }  \ands
g^\ast \a  + \c \in   \sW _{+ ,0} ^{ q \times p }; \label{incluD12} \\
& \d+ g^\ast \b - e_q \in   \sW _{+,0}^{q \times q} \ands
g \d + \b \in  \sW_{-,0}^{p \times q}. \label{incluD34}
\end{align}
Here $ e_p $ and $e_q$ denote the functions identically equal to the identity matrices
$ I_p$  and $I_q$, respectively.  If $g$ has these properties, we refer to $g$ as a
\emph{solution to the twofold EG inverse problem associated with the data set}
$\{  \a, \b, \c,\d \}$.  If a solution exists, then we know from Theorem 1.2  in \cite{KvSch13} that necessarily the following identities hold:
\begin{equation}\label{condD1}
\a^\ast \a - \c^\ast \c =  a_0, \quad \d^\ast \d - \b^\ast \b= d_0,\quad \a^\ast \b = \c^\ast \d.
\end{equation}
Here $a_0$ and $d_0$ are the zero-th Fourier coefficient of $\a$ and $\d$, respectively, and  we
identify the  matrices with $a_0$ and $d_0$ with the matrix functions on $\BT$ that are identically  equal to  $a_0$ and $d_0$, respectively.  Our main problem is to find additional conditions that guarantee the existence of a solution and to obtain explicit formulas for a solution.

The EG inverse problem related to \eqref{incluD12} only  and using $\a$ and $\c$ only  has been treated in \cite{KvSch14}. Here we deal with the  inverse problem   \eqref{incluD12} and  \eqref{incluD34} together,  and for that reason we refer to the problem as a twofold  EG inverse problem. The acronym EG stands for R. Ellis and I. Gohberg,  the authors of \cite{EG92}, where the inverse problem  is solved for the  scalar case, see \cite [Section 4]{EG92}.

Given a data set $\{ \a, \b, \c, \d \} $ and assuming both matrices $a_0$ and $d_0$ are invertible,
 our main theorem (Theorem \ref{thm:mainthmDS1}) gives necessary and sufficient conditions in order that the   twofold EG inverse problem associated with the given data set  has a solution. Furthermore, we show that the solution is unique and we give an explicit formula for the solution in terms of the given data. The results obtained  can be seen as an addition to Chapter 11 in  the Ellis and  Gohberg book  \cite{EG03}.
For some more insight in the role of the matrices $a_0$ and $d_0$ in \eqref{condD1} we refer to
Section \ref{sec:a0d0}.

To understand better the origin of the problem  and to prove our main results we shall restate the twofold EG inverse problem as an  operator problem.
This requires some
further notation and terminology.
 For any positive integer $n$  we denote by $\ell_+^2 ( \BC^{n}) $  and  {$\ell_-^2  ( \BC^{n}) $  the Hilbert spaces
\begin{equation}\label{ell2n} 
\ell^2_{+} ( \BC^{n}) = \Big\{
\begin{bmatrix}x_0\\ x_1\\ x_2\\ \vdots\end{bmatrix}
\mid\sum_{j=0}^\iy \|x_j\| ^2< \infty \Big\},  \ell^2_{-} ( \BC^{n}) = \Big\{
\begin{bmatrix} \vdots \\ x_{-2} \\ x_{-1} \\x_0\\ \end{bmatrix}
\mid\sum_{j=0}^\iy \|x_{-j}\| ^2< \infty \Big\}.   
\end{equation}
We shall also need the corresponding $\ell^1$-spaces which appear when
the superscripts  $2$ in \eqref{ell2n}  are replaced by  $1$.
Since an absolutely summable sequence is square summable, $\ell^1_{\pm} ( \BC^{n})\subset \ell^2_{\pm} ( \BC^{n})$.
In the sequel the one column matrices of the type appearing in \eqref{ell2n} will be denoted by
\[
\begin{bmatrix} x_0 & x_1 & x_2&\cdots  \end{bmatrix}{}^\top\ \mbox{and} \ \begin{bmatrix}\cdots &  x_{-2} & x_{-1}&x_0  \end{bmatrix}{}^\top, \ \mbox{respectively},
\]
with the $\top$-superscript indicating the block transpose. We will also use this notation when the entries are matrices. Finally, let $h$ and $k$ be the linear maps defined by
\begin{align*}
&h =\begin{bmatrix}h_0&h_1&h_2&\cdots\end{bmatrix}{}^\top: \BC^m\to   \ell^1_+(\BC^n), \\[.1cm]
&k =\begin{bmatrix}\cdots&k_{-2}&k_{ -1}&k_{0}\end{bmatrix}{}^\top: \BC^m\to   \ell^1_-(\BC^n).
\end{align*}
With these linear maps we associate the functions    $\sF h$ and $\sF  k$   which are given by
\[
(\sF h)(\l)=\sum_{\nu=0}^\iy \l^\nu h_\nu \quad ( |\l|\leq 1), \quad  (\sF k)(\l)=\sum_{\nu=0}^\iy \l^{-\nu} k_{-\nu} \quad ( |\l|\geq 1).
\]
Since in both cases  the sequence of coefficients are summable in norm, the function $\sF h$ belongs to $\sW_+^{n\ts m}$ and  $\sF k$ belongs to $\sW_-^{n\ts m}$. We shall refer  to $\sF h$ and  $\sF k$ as the \emph{inverse Fourier transforms}  of $h$ and  $k$, respectively.

\medskip
\noindent\textbf{The twofold EG inverse problem as an operator problem.}
Let $\a$, $\b$, $\c$,  $\d$ be the functions appearing in \eqref{Fabcd2}.
With these functions we associate the linear maps:
\begin{align}
&a =\begin{bmatrix}a_0&a_1&a_2&\cdots\end{bmatrix}{}^\top: \BC^p\to   \ell^1_+(\BC^p), \label{alm}\\[.1cm]
&b =\begin{bmatrix}b_0&b_1&b_2&\cdots\end{bmatrix}{}^\top: \BC^q\to   \ell^1_+(\BC^p), \label{blm}\\[.1cm]
&c =\begin{bmatrix}\cdots&c_{-2}&c_{ -1}&c_{0}\end{bmatrix}{}^\top: \BC^p\to   \ell^1_-(\BC^q).\label{clm}\\[.1cm]
&d =\begin{bmatrix}\cdots&d_{-2}&d_{ -1}&d_{0}\end{bmatrix}{}^\top: \BC^q\to   \ell^1_-(\BC^q).\label{dlm}
\end{align}
Here for each $j$ the matrices $a_j$, $b_j$, $c_j$, $d_j$ denote the $j$-th Fourier
coefficients of the functions $\a$, $\b$, $\c$,  $\d$, respectively. Thus   these linear maps are uniquely determined by the functions $\a$, $\b$, $\c$,  $\d$ via the following identities
\begin{equation}\label{defFabcd}
\a=\sF a, \quad  \b=\sF b, \quad  \c=\sF c, \quad \d=\sF d.
\end{equation}
Conversely, if  $a$, $b$, $c$,  $d$ are linear maps as in \eqref{alm} --  \eqref{dlm}, then the functions   $\a$, $\b$, $\c$,  $\d$ defined by \eqref{defFabcd} satisfy the inclusions listed in \eqref{Fabcd2}.

Next,   let $g$ be any function in $\sW_+^{p\ts q}$, say  $g(\z)=\sum_{\nu=0}^\iy  \z^{\nu}g_\nu$, $\z\in \BT$.  With $g$ we associate the Hankel operator $G$ defined by
\begin{equation}\label{defG}
G=\begin{bmatrix}
\cdots&g_2&g_1&g_0\\
\cdots& g_3&g_2&g_1\\
\cdots& g_4&g_3&g_2\\
&\vdots&\vdots&\vdots
\end{bmatrix}: \ell_-^2(\BC^q)\to \ell_+^2(\BC^{p}).
\end{equation}
Using the linear maps $a, b, c, d$ defined by \eqref{alm} -- \eqref{dlm}
it is straightforward to check that
\begin{align}
\eqref{incluD12} \quad  &\Longleftrightarrow  \quad \begin{bmatrix}I&G\\ G^*&I\end{bmatrix}
\begin{bmatrix} a\\ c\end{bmatrix}=\begin{bmatrix} \vp_{+,p}\\0\end{bmatrix}, \label{incluD12a}\\
\eqref{incluD34} \quad  &\Longleftrightarrow  \quad   \begin{bmatrix}I&G\\ G^*&I\end{bmatrix}
\begin{bmatrix} b\\ d\end{bmatrix}=\begin{bmatrix}0\\ \vp_{-,q}\end{bmatrix}. \label{incluD34a}
\end{align}
Here $I$ stands for the identity operator on $\ell_-^2(\BC^q)$ or $\ell_+^2(\BC^{p})$, and
\begin{align*}
\vp_{+,p}&=\begin{bmatrix}I_p&0&0&\cdots\end{bmatrix}{}^\top: \BC^p \to   \ell^1_+(\BC^p), \\[.1cm]
\vp_{-,q}&=\begin{bmatrix}\cdots&0&0&I_q\end{bmatrix}{}^\top: \BC^q\to   \ell^1_-(\BC^q).
\end{align*}
Thus, given the data set $\{\a,\b, \c, \d\}$ and the associate linear maps $a,b,c,d$,  the twofold EG inverse problem is equivalent to  the problem of finding a function  $g \in \sW_+^{p\ts q}$,  $g(\z)=\sum_{\nu=0}^\iy  \z^{\nu}g_\nu$, such that for the Hankel operator $G$ defined by $g$  in \eqref{defG}   the following two identities are satisfied:
\begin{equation}\label{equivprobl}
\begin{bmatrix}I&G\\ G^*&I\end{bmatrix}
\begin{bmatrix} a\\ c\end{bmatrix}=\begin{bmatrix} \vp_{+,p}\\0\end{bmatrix}, \quad
 \begin{bmatrix}I&G\\ G^*&I\end{bmatrix}
\begin{bmatrix} b\\ d\end{bmatrix}=\begin{bmatrix}0\\ \vp_{-,q}\end{bmatrix}.
\end{equation}

In this operator setting the twofold EG inverse problem appears as an infinite dimensional  analogue
of the classical inverse problem for  an $n \ts n$ Hermitian block Toeplitz matrix
$T=[ t_{i-j}]_{i,j=0}^{n}$, where $t_k=t_k^*$, $0\leq k\leq n$,  are $p\ts p$ matrices.
For the latter problem the data consist of two matrix polynomials,
$x(\l)=\sum_{\nu=0}^{n} \l^\nu x_\nu$
and $z(\l)=\sum_{\nu=0}^{n}\l^{-\nu}z_{-\nu}$, with the coefficients being $p\ts p$ matrices, and the  problem is to find $p\ts p$ matrices  $t_0, t_1, \cdots, t_n$ such that
\[
T\begin{bmatrix} x_0 \\ x_1 \\ \vdots \\ x_n \end{bmatrix}=
\begin{bmatrix} I_p \\ 0 \\ \vdots \\ 0 \end{bmatrix} \ands
T\begin{bmatrix}z_n\\  \vdots\\  z_{-1} \\ z_0 \end{bmatrix}=
\begin{bmatrix} 0 \\ \vdots\\ 0 \\ I_p \end{bmatrix}.
\]
The solution of this block Toeplitz matrix inverse problem is due to Gohberg-Heinig \cite{GH74}; see also Theorem 3.4 in \cite{GLe88}.

\medskip
The twofold EG inverse problem is also closely related to an inversion theorem  for  the operator $\om$ defined by the  $2 \ts 2$ operator matrix appearing in \eqref{equivprobl}. Thus
\begin{equation}
\om=\begin{bmatrix}I&G\\ G^*&I\end{bmatrix}, \hspace{.2cm} \mbox{where $G$ is given by \eqref{defG} and $g \in \sW_+^{p\ts q}$}. \label{defM1a}
\end{equation}
In fact, the operator $\om$   is invertible whenever  there exist linear maps $a$, $b$, $c$, $d$ as in \eqref{alm} -- \eqref{dlm} such that the identities in  \eqref{equivprobl} are satisfied. The latter result is given by Theorem 3.1 in \cite{EGL95}, and   without the formula for the inverse of $\om$ it can also be found in Section 11 of the Ellis-Gohberg book \cite{EG03}. The  inversion  theorem also appears  in Section 5 of \cite{FK10} in a more general non-symmetric setting (see  also \cite[Theorem 1.1]{KvSch12}).
We will give a direct proof of  this inversion theorem in Section \ref{sec:prThm3.1}.  The result itself, see Theorem \ref{thm:inversionDS2} in Section~\ref{sec:inversion},  plays
an important role in the  proof  of  our main theorem  (Theorem \ref{thm:mainthmDS1}).

\medskip
\noindent\textbf{Contents.}
The paper consists of ten sections including the present introduction and an appendix.
In Section \ref{sec:LHT} we review a number of standard facts about Laurent, Toeplitz and Hankel operators that are used throughout the paper. In this section we also reformulate the inclusions in \eqref{incluD12} and \eqref{incluD34} in operator language.
In Section \ref{sec:inversion} we state the inversion theorem, and in
Section \ref{sec:solinvprobl} we present the solution of the inverse problem.
In Section \ref{sec:solinvprobl} we also consider the special case when $\det \a$ and $\det \d$ have no zeros in $|\l|\leq 1$ and $|\l|\geq 1$, respectively.
In Section \ref{sec:towardsprf} we prove a number of basic identities that will play a fundamental role in proving Theorem \ref{thm:inversionDS2} in Section \ref{sec:prThm3.1} and
Theorem \ref{thm:mainthmDS1} in Section \ref{sec:prThm4.1}.
In deriving these basic identities we only use that the data $\{\a, \b, \c, \d\}$ satisfy the three conditions in \eqref{condD1} and that the matrices $a_0$ and $d_0$ are invertible.
In Section \ref{sec:Contr} we prove Theorem \ref{thm:posdef}
and in Section~\ref{sec:Polynomial} we
state and prove the solution to the EG inverse problem for the case when
the functions $\a$ and $\b$ are polynomials in $ \lambda$ and the functions
$ \c $ and $\d $ are polynomials in $ \lambda^{-1} $.
In the final section, the appendix, we review (in a somewhat more general setting, allowing the matrices to be non-square) some known results on properties of the matrices $ a_0 $ and $ d_0 $ in the identities in \eqref{condD1},
and present their proofs for the sake of completeness.

\setcounter{equation}{0}
\section{Preliminaries on Hankel and Toeplitz operators} \label{sec:LHT}
We shall need some standard facts  involving  Laurent, Toeplitz,
and Hankel operators  (see, e.g., the first three sections of  \cite[Chapter XXIII]{GGK2}).
Fix   a   $\rho\in \sW^{n\ts m}$, where
$\rho( \z)=\sum_{\nu= - \infty}^\iy r_\nu \z^{\nu}$, $\z\in \BT$.
With $\rho$ we associate an  operator $L_\rho$ which is given by the following $2\ts 2$
operator matrix representation:
\begin{equation}\label{defLrho1}
L_\rho=
\begin{bmatrix} T_{-, \rho}&S_{-,n}^*H_{-, \rho}\\[.2cm] S_{+,n}^*H_{+, \rho}&T_{+, \rho}
 \end{bmatrix}:
\begin{bmatrix}\ell_-^2(\BC^m)\\[.2cm] \ell_+^2(\BC^m)\end{bmatrix} \to
 \begin{bmatrix}\ell_-^2(\BC^{n})\\[.2cm]\ell_+^2(\BC^{n})\end{bmatrix}.
\end{equation}
Here $T_{+,\rho}$ and $T_{-,\rho}$ are the block Toeplitz operators defined by
\begin{align}
T_{+,\rho}&=\begin{bmatrix}
 r_0     &  r_{-1}  &  r_{-2} &\cdots\\
 r_1  & r_0 & r_{-1}&\cdots\\
r_2 & r_1&  r_0 & \\
\vdots      &  \vdots&&\ddots
\end{bmatrix}: \ell_+^2(\BC^m)\to \ell_+^2(\BC^{n}),\label{defTplus}\\
T_{-,\rho}&=\begin{bmatrix}
\ddots &  \vdots      &  \vdots    &     \vdots        \\
\cdots  &r_0 & r_{-1} &   r_{-2}         \\
\cdots        &  r_1 & r_0 & r_{-1}  \\
\cdots &   r_{2}   &  r_{1} & r_0   \end{bmatrix} :\ell_-^2(\BC^m)\to \ell_-^2(\BC^{n}), \label{defTmin}
\end{align}
and  $H_{+,\rho}$ and $H_{-,\rho}$ are the block Hankel operators defined by
\begin{align}
H_{+,\rho}&=\begin{bmatrix}
\cdots& r_2&r_1&r_0\\
\cdots& r_3&r_2&r_1\\
\cdots& r_4&r_3&r_2\\
&\vdots&\vdots&\vdots
\end{bmatrix}: \ell_-^2(\BC^m)\to \ell_+^2(\BC^{n}),\label{defHplus}\\
H_{-,\rho}&=\begin{bmatrix}
\vdots&\vdots&\vdots&\\
r_{-2}&r_{-3}&r_{-4}&\cdots\\
r_{-1}&r_{-2}&r_{-3}&\cdots\\
r_{0}&r_{-1}&r_{-2}&\cdots
\end{bmatrix}: \ell_+^2(\BC^m)\to \ell_-^2(\BC^{n}).\label{defHmin}
\end{align}
Furthermore,  $S_{-,n}$ and $S_{+,n}$ are the block forward shifts  on $\ell_{-}^2(\BC^{n})$ and $\ell_+^2(\BC^{n})$, respectively, that is,
\begin{align}
&S_{-,n}\begin{bmatrix}\cdots &  x_{-2} & x_{-1}&x_0  \end{bmatrix}{}^\top=
\begin{bmatrix}\cdots &  x_{-1} & x_{0}&0  \end{bmatrix}{}^\top, \label{shiftminp}\\
&S_{+,n}\begin{bmatrix} y_0 & y_1 & y_2&\cdots  \end{bmatrix}{}^\top=\begin{bmatrix} 0 & y_0 & y_1&\cdots  \end{bmatrix}{}^\top,  \label{shiftplusq}
\end{align}
and $S_{-,n}^*$ and  $S_{+,n}^*$ are the adjoints of these operators. It follows that
\begin{equation}\label{Hankels1}
S_{-,n}^*H_{-, \rho}=\begin{bmatrix}
\vdots&\vdots&\vdots&\\
r_{-3}&r_{-4}&r_{-5}&\cdots\\
r_{-2}&r_{-3}&r_{-4}&\cdots\\
r_{-1}&r_{-2}&r_{-3}&\cdots
\end{bmatrix}, \quad
S_{+,n}^*H_{+, \rho}=\begin{bmatrix}
\cdots& r_3&r_2&r_1\\
\cdots& r_4&r_3&r_2\\
\cdots& r_5&r_4&r_3\\
&\vdots&\vdots&\vdots
\end{bmatrix}.  
\end{equation}
Moreover, we have
\begin{equation}\label{Hankshift}
S_{-,n}^*H_{-, \rho}=H_{-, \rho}S_{+,m}\ands S_{+,n}^*H_{+, \rho}=H_{+, \rho}S_{-,m}.
\end{equation}
With some  ambiguity  in terminology we call the operator  $L_\rho$ given by \eqref{defLrho1}   the  \emph{Laurent operator} defined by $\rho$.  Since $\rho^*(\zeta)=\sum_{\nu=-\iy}^\iy \z^\nu r_{-\nu}^*$, formulas \eqref{defTplus}, \eqref{defTmin} and  \eqref{defHplus}, \eqref{defHmin} yield the  following identities:
\begin{equation}\label{adjointsTH}
T_{+, \rho}^*=T_{+, \rho^*}, \quad T_{-, \rho}^*=T_{-, \rho^*}, \quad H_{+, \rho}^*=H_{-, \rho^*}.
\end{equation}

In the sequel,   $ \ell $ stands for  the scalar function $ \ell $   given by $ \ell(\z) = \z$ for each $ \z \in \BT $.  Note that $\ell^*(\z)=\ell(\z)^*=\ell(\z)^{-1}$, $ \z  \in \BT $. It follows that
\begin{equation}\label{elrho1}
(\ell\rho)(\z)=\sum_{\nu=-\iy}^\iy \z^\nu r_{\nu-1}\ands  (\ell^* \rho)(\z)=\sum_{\nu=-\iy}^\iy \z^\nu r_{\nu+1} \quad (\z\in \BT).
\end{equation}
But then \eqref{Hankels1} can be rewritten as
\begin{equation}\label{Hankels2}
 S_{-,n}^\ast H_{-, \rho} = H_{-, \ell \rho}\ands S_{+,n}^\ast H_{+, \rho} = H_{+,\ell^\ast \rho}.
\end{equation}

Note that the identities in \eqref{Hankels2} allow us to rewrite  the $2\ts 2$ operator matrix  defining  $L_\rho$ in \eqref{defLrho1} in the following  way:
\begin{equation}\label{defLrho2}
 L_\rho=\begin{bmatrix}T_{-, \rho} &H_{-, \ell \rho} \\   H_{+,\ell^\ast \rho}&T_{+, \rho} \end{bmatrix}.
\end{equation}
Since $\sW=\sW^{1\ts 1}$ is an algebra, we know that  $\rho \in \sW^{ n \ts m}$ and
$\phi \in  \sW^{m\ts k }$ implies  that $\rho \phi \in \sW^{n \ts k}$, and
hence by the theory of Laurent operators we have
\begin{equation}\label{prodrule1}
L_{\rho \phi}=L_{\rho} L_{ \phi}.
\end{equation}
Using the representation  \eqref{defLrho2} for $\rho$, for $\phi$ in place of $\rho$, and for $ \rho\phi$ in place of $\rho$ we see that the product formula \eqref{prodrule1} is equivalent to the following four identities:
\begin{align}
T_{+, \rho \phi }& = T_{+, \rho } T_{+, \phi } + H_{+,\ell^\ast \rho} H_{-, \ell \phi },
\label{Tplusprod} \\
 H_{+, \ell^\ast \rho \phi } &= H_{+,\ell^\ast \rho} T_{-, \phi } +
 T_{+, \rho } H_{+,\ell^\ast \phi },
 \label{Hplusprod} \\
H_{-,\ell \rho \phi } &=  T_{-, \rho} H_{-,\ell \phi } + H_{-, \ell \rho } T_{+, \phi },
\label{Hminprod} \\
 T_{-, \rho \phi }& =  T_{-, \rho} T_{-, \phi } + H_{-,\ell \rho } H_{+,\ell^\ast \phi }.
\label{Tminprod}
\end{align}

Finally, if $r_0$ is an $n \ts n$ matrix, then $ \de_{r_0}$ denotes
the diagonal operator acting on $\ell_-^2(\BC^n)$ or $\ell_+^2(\BC^n)$.
For $ \rho \in \sW^{n \times m } $ one has
\begin{equation}\label{prule}
T_{\pm, \rho r_0 } = T_{\pm , \rho} \Delta_{r_0}  \ands
H_{\pm , \rho r_0 } = H_{\pm , \rho} \Delta_{r_0}.
\end{equation}

\medskip
In the remaining part of this section we deal with the functions $\a$, $\b$, $\c$, $\d$
given by \eqref{Fabcd2}.
The fact  that $ \a \in \sW_+^{p \ts p} $ , $ \b \in \sW_+^{p \ts q} $, $ \c \in \sW_-^{q \ts p} $,  $ \d \in \sW_-^{q \ts q} $ implies  that
\begin{align}
& H_{-,\ell \a} = 0, \quad H_{-,\ell \b} = 0, \quad H_{+,\ell^\ast \c} = 0, \quad
H_{+,\ell^\ast \d} = 0, \label{Hzero1}\\
& H_{+,\ell^\ast \a^\ast } = 0, \quad H_{+,\ell^\ast \b^\ast } = 0,
\quad H_{-, \ell \c^\ast } = 0, \quad H_{-,\ell \d^\ast } = 0. \label{Hzero2}
\end{align}
Note that the identities in \eqref{Hzero2} follow  from those in \eqref{Hzero1} by taking  adjoints.

The next proposition presents some implications of the inclusions in \eqref{incluD12} and \eqref{incluD34} in operator language.

\begin{prop}\label{prop:basicids7}
Let $\a$, $\b$, $\c$, $\d$ be the functions given by \eqref{Fabcd2},  and let $g$
be an arbitrary function in the Wiener space $\sW_+^{p\ts q}$. Then the inclusions in \eqref{incluD12} and \eqref{incluD34} imply the following identities
\begin{align}
& H_{+,g} T_{-, \ell^\ast \c } = - H_{+, \ell^\ast \a } \ands T_{+,\a^\ast } H_{+,g} = - H_{+,\c^\ast}; \label{basicids1a} \\
& T_{+,\ell^\ast \b^\ast } H_{+,g} = - H_{+,\ell^\ast \d^\ast} \ands
H_{+,g} T_{-,\d } = - H_{+,\b}. \label{basicids2a}
\end{align}
More precisely, the first inclusions in \eqref{incluD12} and \eqref{incluD34} imply the first identities in \eqref{basicids1a} and  \eqref{basicids2a},  respectively, and similarly  with second in place of first.
\end{prop}

\bpr We shall only prove that the first inclusion in \eqref{incluD12} implies the first identity in \eqref{basicids1a}. The other implications are proved in a similar way.

First we apply \eqref{Hplusprod} with $ \rho = \ell g $ and $ \phi = \ell^\ast \c$.
Using $\ell$ is scalar we obtain $\rho\phi=\ell g \ell^* \gamma=\ell \ell^* g  \gamma=g\gamma$. This yields
\[
H_{+,g} T_{-,\ell^\ast \c} = H_{+,\ell^\ast g \c} - T_{+, \ell g} H_{+,\ell^\ast \ell^\ast \c}.
\]
Since  $ \ell^\ast \ell^\ast \c \in \sW_{-,0}^{ q \times p } $,  we see that
$ H_{+,\ell^\ast \ell^\ast \c} = 0 $, and thus
\begin{equation}\label{newiden1}
H_{+,g} T_{-,\ell^\ast \c} = H_{+,\ell^\ast g \c}.
\end{equation}
The first inclusion in \eqref{incluD12} tells us that
$ \ell^\ast g \c +  \ell^\ast \a = \ell^\ast ( g \c + \a ) \in \sW_{-,0}^{ p \times p } $, and hence $ H_{+,\ell^\ast g \c} =  -H_{+, \ell^\ast \a } $. Using the latter identity in \eqref{newiden1} we obtain  the first identity in \eqref{basicids1a}.\epr

\setcounter{equation}{0}
\section{The inversion theorem}\label{sec:inversion}
Let $\a$, $\b$, $\c$, $\d$ be the functions  given by \eqref{Fabcd2},  and let $a, b, c, d$ be the associate linear maps given by formulas \eqref{alm}--\eqref{dlm}.
Assume that the two matrices $a_0$ and $d_0$ are invertible.
Using Toeplitz operators and Hankel operators of the type defined in the previous section
we introduce  the following operators:
\begin{align}
M_{11} &= T_{+,\a} \de_{a_{0}^{-1}} T_{+,\a}^* - S_{+,p} T_{+, \b} \de_{d_{0}^{-1}}
T_{+, \b}^*S_{+, p}^* :\ell_+^2(\BC^p) \to \ell_+^2(\BC^p),  \label{defM11} \\[.1cm]
M_{21} &= H_{-,\c} \de_{a_{0}^{-1}} T_{+,\a}^* - S_{-,q}^* H_{-, \d} \de_{d_{0}^{-1}}
T_{+, \b}^* S_{+, p}^*  : \ell_+^2(\BC^p) \to \ell_-^2(\BC^q),\label{defM21} \\[.1cm]
M_{12}&= H_{ +, \b} \de_{d_0^{-1}} T_{-, \d}^* - S_{+,p}^*H_{+, \a} \de_{a_0^{-1}}
T_{-,\c}^* S_{-,q}^* : \ell_-^2(\BC^q) \to \ell_+^2(\BC^p),
\label{defM12} \\[.1cm]
M_{22}&= T_{-,\d}\de_{d_{0}^{-1}} T_{-, \d}^* -  S_{-,q}T_{-, \c}\de_{a_{0}^{-1}}T_{-,\c}^*S_{-, q}^* :\ell_-^2(\BC^q)\to \ell_-^2(\BC^q). \label{defM22}
\end{align}
Notice that the operators $ M_{ij} $, $ 1 \leq i,j \leq 2 $, are uniquely determined by the data.
If $a_0$ and $d_0$ are selfadjoint, then formulas \eqref{defM11} and \eqref{defM22} show  that $ M_{11}^\ast = M_{11} $ and $ M_{22}^\ast = M_{22} $. Later (see Lemma \ref{lem:altMij}) we shall prove that under certain additional conditions $ M_{12}^\ast = M_{21}$.

We are now ready to state the inversion theorem.

\begin{thm}\label{thm:inversionDS2}
Let $\a$, $\b$, $\c$, $\d$ be the functions  given by \eqref{Fabcd2},  and let $a, b, c, d$ be the associate linear maps given by formulas \eqref{alm}--\eqref{dlm} with both matrices $a_0$ and $d_0$ invertible.
Assume $g\in\sW_+^{p\ts q} $ is a solution to the twofold EG inverse problem associated with the data set $\{\a, \b, \c, \d\}$.
Then $g$ is the only solution and   the  operator  $\om $ given by the $2 \times 2 $ operator matrix
\begin{equation}\label{defM1b}
\om =\begin{bmatrix} I &H_{+, g}\\ H_{-, g^*}&I\end{bmatrix}:
\begin{bmatrix} \ell_+^2(\BC^p)\\ \ell_-^2(\BC^q) \end{bmatrix}\to
\begin{bmatrix} \ell_+^2(\BC^p)\\ \ell_-^2(\BC^q)  \end{bmatrix}
\end{equation}
is invertible and its inverse is given by
\begin{equation}\label{formO-1}
\om^{-1} = M  =\begin{bmatrix} M_{11}&M_{12}\\ M_{21}&M_{22}\end{bmatrix},
\end{equation}
where $M_{11}$, $M_{21}$, $M_{12}$, $M_{22}$ are the operators given by \eqref{defM11} -- \eqref{defM22}.

Conversely, if $g\in\sW_+^{p\ts q} $, and the operator $ \om $ given by \eqref{defM1b} is invertible, then there exists a unique data set  $\{\a, \b, \c, \d\}$
such that $g$ is a solution to the twofold EG inverse problem associated with this data set.
\end{thm}

\begin{rem}\label{rem:inversionthm}
\textup{As is mentioned in the introduction, the above theorem is known.
In Section \ref{sec:prThm3.1} we shall give a direct proof based on the analysis of the operators
$M_{ij}$, $1\leq i,j\leq 2$, given in Section \ref{sec:towardsprf}. \  }
\end{rem}

For later purposes we mention that the operators $M_{ij}$, $1\leq i,j\leq 2$, are also given by
\begin{align}
M_{11} &= T_{+,\a} \de_{a_{0}^{-1}} T_{+,\a}^\ast -  T_{+,\ell \b} \de_{d_{0}^{-1}}
T_{+,\ell \b}^\ast : \ell_+^2(\BC^p) \to \ell_+^2(\BC^p),  \label{def2M11} \\[.1cm]
M_{21} &= H_{-,\c} \de_{a_{0}^{-1}} T_{+,\a}^\ast - H_{-,\ell \d} \de_{d_{0}^{-1}}
T_{+,\ell \b}^\ast  : \ell_+^2(\BC^p) \to \ell_-^2(\BC^q), \label{def2M21} \\[.1cm]
M_{12}&= H_{+, \b} \de_{d_0^{-1}} T_{-, \d}^\ast - H_{+,\ell^\ast \a} \de_{a_0^{-1}}
T_{-,\ell^\ast \c}^\ast : \ell_-^2(\BC^q) \to \ell_+^2(\BC^p),
\label{def2M12} \\[.1cm]
M_{22}&= T_{-,\d} \de_{d_{0}^{-1}} T_{-, \d}^\ast -  T_{-,\ell^\ast \c}
\de_{a_{0}^{-1}} T_{-,\ell^\ast \c}^\ast :\ell_-^2(\BC^q)\to \ell_-^2(\BC^q). \label{def2M22}
\end{align}
Here $\ell$ is the  scalar function defined in the paragraph directly after \eqref{adjointsTH}.

To derive the above formulas note  that
\begin{align*}
&S_{+,p} T_{+, \b} =S_{+,p} \begin{bmatrix}
b_0     & 0 & 0 &\cdots\\
b_1  & b_0 & 0&\cdots\\
b_2 & b_1&  b_0 & \\
\vdots      &  \vdots&&\ddots
\end{bmatrix}=\begin{bmatrix}
0     & 0 & 0 &\cdots\\
b_0  & 0 & 0&\cdots\\
b_1 & b_0&  0 & \\
\vdots      &  \vdots&&\ddots
\end{bmatrix}=T_{+,\ell \b},\\
&S_{-,q} T_{-, \c}=S_{-,q}\begin{bmatrix}
\ddots &  \vdots      &  \vdots    &     \vdots        \\
\cdots  &c_0 & c_{-1} &   c_{-2}         \\
\cdots   &  0 & c_0 & c_{-1}  \\
\cdots &  0   & 0 & c_0   \end{bmatrix} = \begin{bmatrix}
\ddots &  \vdots      &  \vdots    &     \vdots        \\
\cdots  &0& c_{0} &   c_{-1}         \\
\cdots   &  0 & 0 & c_{0}  \\
\cdots &  0   & 0 & 0  \end{bmatrix}=T_{-,\ell^* \c}.
\end{align*}
Furthermore, from the identities in  \eqref{Hankels2} we know  that
\[
S_{-,q}^* H_{-, \d}=H_{-, \ell\d} \ands S_{+,p}^* H_{+, \a}= H_{+, \ell^*\a}.
\]
Using the above identities in \eqref{defM11}--\eqref{defM22} we obtain \eqref{def2M11}--\eqref{def2M22}.

\begin{rem}\label{a0d0inv}
\textup{In Theorem  \ref{thm:inversionDS2} the condition that both  $a_0$ and $d_0$ are invertible can be replaced by the weaker condition that at least one of the two is invertible. On the other hand, in that case the assumption $g\in\sW_+^{p\ts q} $ is a solution to the twofold EG inverse problem  implies that both are invertible; see, e.g.,
\cite[Section 3 and 4]{EGL95}, \cite[Proposition 5.2]{FK10}  or
\cite[Theorem 1.1]{KvSch12},
and see Section \ref{sec:a0d0} for further details.}
\end{rem}

\setcounter{equation}{0}
\section{The solution of the inverse problem}\label{sec:solinvprobl}

In this section we present our solution to the twofold EG inverse problem, as well as a characterization of the case where a solution exists with a strictly contractive Hankel operator.

\begin{thm}\label{thm:mainthmDS1}
Let $\a$, $\b$, $\c$, $\d$ be the functions  given by \eqref{Fabcd2},
and let $a, b, c, d$ be the associate linear maps given by formulas \eqref{alm}--\eqref{dlm}
with both matrices $a_0$ and $d_0$ invertible.
Then  the twofold  EG inverse problem associated with the data  set
$ \{ \a, \b, \c, \d \} $
has a solution if and only the following conditions are satisfied:
\begin{itemize}
\item[\textup{(D1)}] the identities in \eqref{condD1} hold true;
\item[\textup{(D2)}]
the operators $ M_{11} $ and $ M_{22} $  defined by  \eqref{defM11}
and \eqref{defM22} are one-to-one.
\end{itemize}
Furthermore, in that case  $ M_{11} $ and $ M_{22} $ are invertible, the solution is unique and  the unique solution $g$ and its adjoint are  given by
\begin{equation} \label{unique1}
g = -    \sF  (M_{11}^{-1}b) \ands    g^\ast  = - \sF (M_{22}^{-1}c).
\end{equation}
\end{thm}

\noindent
\begin{rem} \textup{We shall see  that $M_{11}^{-1}b$ is a linear map from $\BC^q$ to  $\ell_+^1 (\BC^p)$ and $M_{22}^{-1} c$ is a linear map from $\BC^p$ to $\ell_-^1(\BC^p)$. Thus the inverse Fourier transforms in \eqref{unique1} are well-defined.}
\end{rem}

\begin{rem} \textup{The necessity of condition (D1) we know from \cite[Theorem 5]{KvSch13}. Furthermore, we can use Theorem \ref{thm:inversionDS2} to prove the necessity of condition (D2), the uniqueness of the solution and  a formula for
the solution.
Indeed, as we shall see in Section \ref{sec:prThm3.1},  formula
\eqref{formO-1} implies that $ H_{+,g} = -M_{11}^{-1} M_{12} $.
New in the above theorem  are  the sufficiency of conditions (D1) and (D2),  and the
formulas for the solution given in \eqref{unique1}}.
\end{rem}

The next theorem is a generalisation of  Theorem 5.1 in \cite{EGL95} which deals with the inverse problem for the onefold case; see also \cite[Theorem 3.1]{KvSch14}.

\begin{thm}\label{thm:posdef}
Let $\a$, $\b$, $\c$, $\d$ be the functions  given by \eqref{Fabcd2}, and let $a, b, c, d$ be the associate linear maps given by formulas \eqref{alm}--\eqref{dlm}.
Then the twofold  EG inverse problem associated with the data  set
$ \{ \a, \b, \c, \d \} $ has a solution $g$ with the additional property that $H_{+,g}$  is a strict contraction if and only if the following  three conditions are satisfied:
\begin{itemize}
\item[\textup{(i)}]  $a_0$ and $d_0$ are positive definite,
\item[\textup{(ii)}]  $\a^\ast \a - \c^\ast \c =  a_0 $, \ $\d^\ast \d - \b^\ast \b = d_0$,\  $\a^\ast \b = \c^\ast \d $,
\item[\textup{(iii)}]  $\det \a$ and $\det \d$ have no zeros in $|\l|\leq 1$ and $|\l|\geq 1$, respectively.
\end{itemize}
\end{thm}

\setcounter{equation}{0}
\section{Towards the proofs of the theorems}\label{sec:towardsprf}
Throughout this section $\a$, $\b$, $\c$, $\d$ are the functions  given by \eqref{Fabcd2},
and  $a, b, c, d$ are the associate linear maps given by formulas \eqref{alm}--\eqref{dlm}.
Moreover it will be assumed that the matrices $a_0$ and $d_0$ are invertible.
In the four lemmas presented in this section we also assume
that the identities in  \eqref{condD1} are satisfied,  that is,
\begin{equation}\label{condD1a}
\a^\ast \a - \c^\ast \c =  a_0, \quad \d^\ast \d - \b^\ast \b= d_0,\quad \a^\ast \b = \c^\ast \d.
\end{equation}
The first two identities in \eqref{condD1a} imply that the matrices $a_0$ and $ d_0 $ are selfadjoint.
Hence  formulas \eqref{defM11} and \eqref{defM22} give that $ M_{11}^\ast = M_{11} $ and
$ M_{22}^\ast = M_{22}$  (cf., the comment after \eqref{defM22}).

Together the identities in \eqref{condD1a} are equivalent to
\begin{equation}\label{C1-3matrix}
\begin{bmatrix} \a^\ast & \c^\ast  \\ \b^\ast & \d^\ast \end{bmatrix}
\begin{bmatrix} I_p & 0 \\ 0 & - I_q \end{bmatrix}
\begin{bmatrix} \a  & \b \\ \c & \d \end{bmatrix}=
\begin{bmatrix} a_0 & 0 \\ 0 & - d_0 \end{bmatrix}.
\end{equation}
Since all entries are matrices of functions, it follows that
\begin{equation}\label{C1-3matrix*}
\begin{bmatrix} \a  & \b  \\ \c &  \d  \end{bmatrix}
\begin{bmatrix} a_0^{-1}  & 0 \\ 0 & - d_0^{-1} \end{bmatrix}
\begin{bmatrix} \a^\ast & \c^\ast  \\ \b^\ast & \d^\ast \end{bmatrix}
= \begin{bmatrix} I_p & 0 \\ 0 & - I_q \end{bmatrix},
\end{equation}
which in turn is  equivalent to the following three identities
\begin{align}
&\a a_0^{-1} \a^\ast - \b d_0^{-1} \b^\ast  = I_p, \quad
\d d_0^{-1} \d^\ast - \c a_0^{-1} \c^\ast = I_q , \label{conD3a}\\
&\hspace{2.5cm}  \a a_0^{-1} \c^\ast = \b d_0^{-1} \d^\ast.\label{conD3b}
\end{align}
By a similar argument, it follows that the identities in \eqref{conD3a} and  \eqref{conD3b} imply those in \eqref{condD1a}, hence we conclude that the three identities in \eqref{conD3a} and  \eqref{conD3b} are equivalent to the three identities in \eqref{condD1a}
(provided, as is assumed throughout this section, $a_0$ and $d_0$ are invertible).

The following lemma is an addition to the final part of Section \ref{sec:LHT}.

\begin{lem}\label{basicids5}
Assume that  condition \eqref {condD1a} is satisfied. Then
\begin{align}
&\hspace{1cm}\begin{bmatrix}  T_{+,\a^\ast} \\ T_{+, \ell^\ast  \b^\ast}\end{bmatrix}
\begin{bmatrix}  H_{+, \ell^\ast \a} & H_{+,\b} \end{bmatrix} =
\begin{bmatrix}  H_{+,\c^\ast} \\ H_{+, \ell^\ast \d^\ast}\end{bmatrix}
\begin{bmatrix} T_{-, \ell^\ast  \c} &  T_{-,\d} \end{bmatrix},\label{matrixTHHT}\\[.2cm]
&\begin{bmatrix} T_{+,\a^\ast} \\ T_{+,\ell^* \b^*}  \end{bmatrix} S_{+,p}^\ast
 \begin{bmatrix} H_{+, \ell^\ast \a}& H_{+,\b} \end{bmatrix}  =
\begin{bmatrix} H_{+,\c}^\ast \\ H_{+, \ell^\ast \d^\ast} \end{bmatrix} S_{-,q}
\begin{bmatrix} T_{-,\d} & T_{-, \ell^\ast \c } \end{bmatrix}. \label{rhs3}
\end{align}
\end{lem}
\bpr
In order to  prove \eqref{matrixTHHT} we have  to check the following four identities:
\begin{align}
 T_{+, \a^\ast} H_{+, {\ell^\ast \a }} = H_{+, \c^\ast } T_{-, \ell^\ast \c } &\ands
T_{+,\a^\ast } H_{+,\b} = H_{+, \c^\ast}  T_{-,\d},  \label{equal53a}\\
 T_{+, \ell^\ast  \b^\ast}H_{+, \ell^\ast \a}=H_{+, \ell^\ast \d^\ast}T_{-, \ell^\ast  \c} &\ands T_{+, \ell^\ast  \b^\ast}H_{+,\b}=H_{+, \ell^\ast \d^\ast}T_{-,\d}.\label{equal53b}
\end{align}

\smallskip\noindent
\textsc{Step 1.}
We  prove   the second identity in \eqref{equal53a}.   First we  apply \eqref{Hplusprod}  with $ \rho =\a^\ast $ and $ \phi = \ell \b$, and  we use  the first identity in  \eqref{Hzero2}. We  obtain
\[
H_{+, \a^\ast \b } = H_{+, \ell^\ast \a^\ast \ell \b } = H_{+, \ell^\ast \a^\ast } T_{-,\ell \b} +
T_{+,\a^\ast } H_{+,\b} = T_{+,\a^\ast } H_{+,\b}.
\]
On the other hand,  by applying \eqref{Hplusprod}  with  $ \rho =\ell \c^\ast $ and $ \phi = \d $,   and using the last equality in \eqref{Hzero1},  we get
\[
H_{+, \c^\ast \d } = H_{+, \ell^\ast  \ell \c^\ast \d }= H_{+, \c^\ast} T_{-,\d} +
T_{+, \ell \c^\ast } H_{+, \ell^\ast \d} =  H_{+, \c^\ast}  T_{-,\d}.
\]
Now notice that  the third  identity in \eqref{condD1a} implies  that $ H_{+, \a^\ast \b } = H_{+, \c^\ast \d } $.
It follows that   $T_{+,\a^\ast } H_{+,\b} = H_{+, \c^\ast}  T_{-,\d} $ as desired.

\smallskip\noindent
\textsc{Step 2.}  We  prove   the first identity in \eqref{equal53a}.  We first apply   \eqref{Hplusprod}  with $ \rho =\a^\ast $ and $ \phi = \a$, and  we use again the first identity in  \eqref{Hzero2}. This yields
\[
H_{+, \ell^\ast \a^\ast \a } = H_{+, \ell^\ast \a^\ast } T_{-,\a } +
T_{+, \a^\ast} H_{+, \ell^\ast \a } = T_{+, \a^\ast} H_{+, \ell^\ast \a } .
\]
On the other hand, using the first identity in  \eqref{condD1a}, we see that
$ \ell^\ast \a^\ast \a = \ell^\ast \c^\ast \c + \ell^\ast a_0 $. Since the Hankel operator
$ H_{+, \ell^\ast a_0}$ is zero,  we obtain
\[
H_{+,\ell^\ast \a^\ast \a } = H_{+, \ell^\ast \c^\ast \c } .
\]
Next we apply \eqref{Hplusprod}  with $ \rho =\ell \c^\ast $ and $ \phi = \ell^*\c$, which yields
\[
H_{+, \ell^\ast \c^\ast \c }  = H_{+, \ell^\ast \ell \c^\ast \ell^\ast \c }=
H_{+, \c^\ast } T_{-, \ell^\ast \c } + T_{+, \ell \c^\ast} H_{+, \ell^\ast \ell^\ast \c }.
\]
Note that the Hankel operator  $H_{+, \ell^\ast \ell^\ast \c }$ is zero. We conclude that
\[
T_{+, \a^\ast} H_{+, \ell^\ast \a } =H_{+,\ell^\ast \a^\ast \a } = H_{+, \ell^\ast \c^\ast \c } =H_{+, \c^\ast } T_{-, \ell^\ast \c }.
\]
We proved that $T_{+, \a^\ast} H_{+,\ell^\ast \a } = H_{+, \c^\ast } T_{-, \ell^\ast \c}$.

\smallskip\noindent
\textsc{Step 3.} We  prove  the first identity in \eqref{equal53b}.
First we apply \eqref{Hplusprod} with $ \rho = \ell^\ast \b^\ast $ and $ \phi = \a$, and we use that $ H_{+, \ell^\ast \ell^\ast \b^\ast } = 0 $.
This yields
\[
H_{+, \ell^\ast \ell^\ast \b^\ast \a } =
H_{+, \ell^\ast \ell^\ast \b^\ast } T_{-,\a}+
T_{+, \ell^\ast \b^\ast } H_{+,\ell^\ast \a }
= T_{ +, \ell^\ast \b^\ast } H_{+,\ell^\ast \a } .
\]
Next we apply \eqref{Hplusprod} with $ \rho = \d^\ast $ and $ \phi = \ell^\ast \c  $, and we use that $ H_{+, \ell^\ast \ell^\ast \c } = 0 $.
This yields
\[
H_{+, \ell^\ast  \d^\ast \ell^\ast \c } =
H_{+, \ell^\ast \d^\ast } T_{-,\ell^\ast \c} +
T_{+, \d^\ast } H_{+,\ell^\ast \ell^\ast \c }
=H_{+, \ell^\ast \d^\ast } T_{-,\ell^\ast \c} .
\]
Notice that the third identity in \eqref{condD1a} implies that
$ H_{+, \ell^\ast \ell^\ast \b^\ast \a }  = H_{+, \ell^\ast \ell^\ast \d^\ast \c }  $.
Thus we proved   the first identity in \eqref{equal53b}.

\smallskip\noindent
\textsc{Step 4.}
We  prove the second  identity in \eqref{equal53b}.
To do this we  apply \eqref{Hplusprod} with $ \rho = \ell^\ast \b^\ast $ and $ \phi = \ell \b $,
and we use that $ H_{+, \ell^\ast \ell^\ast \b^\ast } = 0 $.
This yields
\[
H_{+, \ell^\ast  \b^\ast \b } = H_{+, \ell^\ast\ell^\ast\beta^\ast\ell \beta }
= H_{+, \ell^\ast \ell^\ast \b^\ast } T_{-,\ell \b}+ T_{+, \ell^\ast \b^\ast } H_{+,\b }
=T_{+, \ell^\ast \b^\ast } H_{+,\b }.
\]
Next we apply \eqref{Hplusprod}  with $ \rho = \d^\ast $ and $ \phi = \d $.
This yields
\[
H_{+, \ell^\ast \d^\ast \d }
= H_{+, \ell^\ast \d^\ast } T_{-,\d }+ T_{+, \d^\ast } H_{+, \ell^\ast \d }
=H_{+, \ell^\ast \d^\ast } T_{-,\d} ,
\]
where the last equality follows from $ H_{+, \ell^\ast \d } = 0 $.
To complete the proof of this step we use the second identity in \eqref{equal53b} to
show that
\[
H_{+,\ell^\ast \d^\ast \d } - H_{+,\ell^\ast \b^\ast \b } = H_{+,\ell^\ast d_0 } =0.
\]

\smallskip\noindent
\textsc{Step 5.}
It remains to prove \eqref{rhs3}. Since $H_{+, \ell^\ast \a}$ and $ H_{+,\b}$ are Hankel operators,  and the operators $T_{-,\d}$ and $T_{-, \ell^\ast \c } $ are Toeplitz operators,  the following intertwining relations hold:
\begin{align*}
&S_{+,p}^\ast
 \begin{bmatrix} H_{+, \ell^\ast \a}& H_{+,\b} \end{bmatrix}  =
 \begin{bmatrix} H_{+, \ell^\ast \a}& H_{+,\b} \end{bmatrix}
 \begin{bmatrix} S_{-q}&0\\0&S_{-p}\end{bmatrix}\\
& S_{-,q} \begin{bmatrix} T_{-,\d} & T_{-, \ell^\ast \c } \end{bmatrix}=
  \begin{bmatrix} T_{-,\d} & T_{-, \ell^\ast \c } \end{bmatrix}\begin{bmatrix} S_{-q}&0\\0&S_{-p}\end{bmatrix}.
\end{align*}
But then using \eqref{matrixTHHT} we obtain:
\begin{align*}
& \begin{bmatrix}T_{+,\a^\ast} \\ T_{+,\ell^* \b^*}  \end{bmatrix} S_{+,p}^\ast
 \begin{bmatrix} H_{+, \ell^\ast \a}& H_{+,\b} \end{bmatrix}  =\\
&\hspace{1cm}= \begin{bmatrix}T_{+,\a^\ast} \\ T_{+,\ell^* \b^*}\end{bmatrix}  \begin{bmatrix} H_{+, \ell^\ast \a}& H_{+,\b} \end{bmatrix}
 \begin{bmatrix} S_{-q}&0\\0&S_{-p}\end{bmatrix}\\
&\hspace{1cm}=\begin{bmatrix}  H_{+,\c^\ast} \\ H_{+, \ell^\ast \d^\ast}\end{bmatrix}
\begin{bmatrix} T_{-, \ell^\ast  \c} &  T_{-,\d} \end{bmatrix}  \begin{bmatrix} S_{-q}&0\\0&S_{-p}\end{bmatrix}\\
&\hspace{1cm}=
\begin{bmatrix}  H_{+,\c^\ast} \\ H_{+, \ell^\ast \d^\ast}\end{bmatrix} S_{-q}
\begin{bmatrix} T_{-, \ell^\ast  \c} &  T_{-,\d} \end{bmatrix},
\end{align*}
and  \eqref{rhs3} is proved.
\epr

\medskip
By taking adjoints and using the identities in \eqref{adjointsTH} it is straightforward  to  see that  \eqref{matrixTHHT} yields the following identity:
\begin{equation}\label{matrixTHHT2}
\begin{bmatrix}  T_{-,\d^\ast} \\ T_{-,\ell \c^\ast}\end{bmatrix}
\begin{bmatrix} H_{-,\c} & H_{-,\ell \d} \end{bmatrix} =
\begin{bmatrix}  H_{-,\b^\ast} \\ H_{-,\ell \a^\ast}\end{bmatrix}
\begin{bmatrix} T_{+,\a} & T_{+,\ell \b} \end{bmatrix}.
\end{equation}
Furthermore, note that the identity \eqref{rhs3} remains true if  $S_{+,p}^*$ and $S_{-,q}$ are
replaced by  $(S_{+,p}^*)^n$ and $(S_{-,q})^n$, respectively, where $n$ is any nonnegative integer.

\medskip
The following three lemmas contain basic identities for the operators
$M_{ij} $, defined by \eqref{defM11}--\eqref{defM22}.
These identities will play an essential role in the proofs of Theorems \ref{thm:inversionDS2} and
\ref{thm:mainthmDS1}.
In fact, they will allow us to reduce the proofs of those two theorems to a matter of direct checking.

\begin{lem}\label{lem:unitsM}
Assume that condition \eqref {condD1a}   is satisfied. Then the following identities hold true:
\begin{equation} \label{unitsM}
M_{11} \vp_{+,p} = a, \quad  M_{12} \vp_{-,q} = b, \quad M_{21} \vp_{+,p} = c,
\quad  M_{22} \vp_{-,q} = d.
\end{equation}
\end{lem}

\bpr
From the  block  matrix representation of $ T_{+,\ell \b} $ one sees that the first column  of  $ T_{+,\ell \b}^\ast $ consists of zero entries only,  and thus  $ T_{+,\ell \b}^\ast \vp_{+,p} =0 $.   Using this fact together with $ a_0^\ast =  a_0 $ in   \eqref{def2M11}  yields
\begin{align*}
M_{11} \vp_{+,p}
&=  T_{+,\a} \de_{a_{0}^{-1}} T_{+,\a}^\ast \vp_{+,p} =
T_{+,\a} \de_{ a_0^{-1} } \vp_{+,p} a_0^\ast =
T_{+,\a} \vp_{+,p} a_0^{-1} a_0 \\
&=T_{+,\a} \vp_{+,p} = a.
\end{align*}
This proves the first identity in \eqref{unitsM}. Again using $ T_{+,\ell \b}^\ast e_{+,p} =0  $ and \eqref{def2M21}  one obtains
\begin{align*}
M_{21} \vp_{+,p}
&= H_{-,\c} \de_{a_{0}^{-1}} T_{+,\a}^\ast \vp_{+,p} =
 H_{-,\c} \de_{a_{0}^{-1}}  \vp_{+,p}a_0^\ast =
 H_{-,\c} \vp_{+,p}a_{0}^{-1} a_0 \\
&=H_{-,\c} \vp_{+,p} = c,
\end{align*}
which proves the third identity in \eqref{unitsM}. The other two identities are proved in a similar way.
\epr

\begin{lem}\label{lem:altMij}
Assume that  condition \eqref {condD1a}  is satisfied.  Then
\begin{align}
& M_{11} = I_{\ell_+^2 (\BC^p) } - H_{+,\ell^\ast \a} \Delta_{a_0^{-1}} H_{+,\ell^\ast \a}^\ast +
H_{+,\b} \Delta_{d_0^{-1}} H_{+,\b }^\ast   ,
\label{formulaM11} \\
& M_{21} =  T_{-,\d} \Delta_{d_0^{-1}} H_{+,\b }^\ast -
T_{-,\ell^\ast \c} \Delta_{a_0^{-1}} H_{+, \ell^\ast \a}^\ast
 , \label{formulaM21}\\
& M_{12} =  T_{+,\a} \Delta_{a_0^{-1}} H_{-,\c}^\ast -
T_{+,\ell \b} \Delta_{d_0^{-1}}  H_{-,\ell \d}^\ast
  , \label{formulaM12} \\
& M_{22} = I_{\ell_-^2 (\BC^q)} - H_{-,\ell \d} \Delta_{d_0^{-1}} H_{-,\ell \d}^\ast +
H_{-,\c} \Delta_{a_0^{-1}} H_{-,\c} ^\ast .
\label{formulaM22}
\end{align}
In particular $ M_{12}^\ast = M_{21} $.   Hence $M$ in \eqref{formO-1} is selfadjoint.
\end{lem}

\bpr
First  notice that condition \eqref{condD1a} yields the first identity in \eqref{conD3a},  and therefore
\begin{equation}\label{iden1a}
T_{+,\a a_0^{-1} \a^\ast } - T_{+,  \b d_0^{-1} (\b)^\ast }-
I_{\ell_+^2 (\BC^p) }=T_{+, \a a_0^{-1} \a^\ast - \b d_0^{-1} \b^\ast - \vp_{+,p} } = 0 .
\end{equation}
Since $ \ell \ell^\ast $ is identically equal to 1, we have
$\ell \b d_0^{-1} (\ell \b)^\ast = \b \ell d_0^{-1} \ell^\ast \b^\ast =
\b d_0^{-1} \b^\ast$. By using this in \eqref{iden1a} we obtain
\begin{equation*}
T_{+,\a a_0^{-1} \a^\ast } - T_{+,\ell  \b d_0^{-1} (\ell \b)^\ast } =
I_{\ell_+^2 (\BC^p) } .
\end{equation*}
By applying the product rule \eqref{Tplusprod} and the identities in \eqref{prule} one sees that
\begin{align*}
T_{+,\a } \Delta_{a_0^{-1}} T_{+,\a }^\ast & -
 T_{+,\ell \b} \Delta_{d_0^{-1}} T_{+,\ell \b}^\ast = \\
&=  I_{\ell_+^2 (\BC^p) } - H_{+, \ell^\ast \a} \Delta_{a_0^{-1}} H_{+, \ell^\ast \a}^\ast
+ H_{+, \b} \Delta_{d_0^{-1}} H_{+,\b}^\ast.
\end{align*}
We conclude that the operator $ M_{11}$  defined  by \eqref{def2M11} is also given  by \eqref{formulaM11}.

In a similar way one shows that the second identity in \eqref{conD3a}
yields the identity \eqref{formulaM22}.

Next we apply \eqref{Hplusprod} with $ \rho = \a a_0^{-1} $ and $ \phi = \ell \c^\ast $. This yields
\begin{align}
H_{+, \a a_0^{-1} \c^\ast }& = H_{+,  \ell^\ast \a a_0^{-1} \ell \c^\ast } \nonumber \\
&= T_{+,\a} \Delta_{a_0^{-1} } H_{+,\c^\ast} +
H_{+,  \ell^\ast \a } \Delta_{a_0^{-1}} T_{-,\ell \c^\ast} \nonumber \\
&= T_{+,\a} \Delta_{a_0^{-1} } H_{-,\c}^\ast +
H_{+,  \ell^\ast \a } \Delta_{a_0^{-1}} T_{-, \ell^\ast \c}^\ast. \label{rhs1}
\end{align}
Similarly,  by applying  \eqref{Hplusprod} with $ \rho = \ell \b d_0^{-1} $
and $ \phi = \d^\ast $, we  obtain
\begin{align}
H_{+, \b d_0^{-1} \d^\ast }& = H_{+,  \ell^\ast  \ell \b d_0^{-1} d^\ast } \nonumber \\
&= T_{+,\ell \b} \Delta_{d_0^{-1} } H_{+, \ell^\ast \d^\ast} +
H_{+, \b } \Delta_{d_0^{-1}} T_{-, \d^\ast} \nonumber \\
&= T_{+,\ell \b} \Delta_{d_0^{-1} } H_{-,\ell \d}^\ast +
H_{+, \b} \Delta_{d_0^{-1}} T_{-, \d}^\ast. \label{rhs2}
\end{align}
Now note that  \eqref{conD3b} implies  that   $ H_{+, \a a_0^{-1} \c^\ast } = H_{+, \b d_0^{-1} \d^\ast } $.  But then  equalities   \eqref{rhs1} and \eqref{rhs2}   show   that the right hand sides of \eqref{def2M12} and \eqref{formulaM12} are equal.

From \eqref{def2M21} and \eqref{formulaM12} we obtain $M_{21}^*=M_{12}$. Using this fact along with formula \eqref{def2M12} for $M_{12}$ we see that $M_{21}=M_{12}^*$ is given by \eqref{formulaM21}. \
We already know (see the first paragraph of the present section)  that $M_{11}$ and $M_{22}$ are selfadjoint.  Hence $M$ in \eqref{formO-1} is selfadjoint.
\epr

\medskip

\begin{lem}\label{lem:inverseR}
Assume that condition \eqref{condD1a} is satisfied.
Let the operators $ M_{ij} $, $ i,j=1,2$,  be given by \eqref{defM11}--\eqref{defM22}. Then
\begin{equation}\label{equality1}
\begin{bmatrix} M_{11} & M_{12} \\ M_{21} &  M_{22}  \end{bmatrix}
\begin{bmatrix} I_{\ell_+^2 (\BC^p) } & 0 \\ 0 &  - I_{\ell_-^2 (\BC^q) } \end{bmatrix}
\begin{bmatrix} M_{11} & M_{12} \\ M_{21} &  M_{22}  \end{bmatrix}  =
\begin{bmatrix} M_{11} & 0 \\ 0 &  -M_{22}  \end{bmatrix}.
\end{equation}
In particular,
\begin{equation}\label{invform1}
M = \begin{bmatrix} M_{11} & M_{12} \\ M_{21} &  M_{22}  \end{bmatrix} :
\begin{bmatrix} \ell_+^2 (\BC^p) \\ \ell_-^2 (\BC^q) \end{bmatrix}  \to
\begin{bmatrix} \ell_+^2 (\BC^p) \\ \ell_-^2 (\BC^q) \end{bmatrix}
\end{equation}
is invertible if and only if $ M_{11} $ and $ M_{22} $ are invertible, and in that case
\begin{equation}\label{invform2}
M^{-1}  =
\begin{bmatrix} I_{\ell_+^2(\BC^p) } & -M_{12} M_{22}^{-1} \\ -M_{21} M_{11}^{-1} &
I_{\ell_-^2(\BC^q) } \end{bmatrix}
=
\begin{bmatrix} I_{\ell_+^2(\BC^p) } & -M_{11}^{-1} M_{12} \\ -M_{22}^{-1} M_{21} &
I_{\ell_-^2(\BC^q)}  \end{bmatrix}.
\end{equation}
Finally,
\begin{equation}\label{intertwiM}
 M_{11} S_{+,p}^\ast M_{12}= M_{12} S_{-,q} M_{22}.
\end{equation}
\end{lem}

\bpr
To check \eqref{equality1} we will prove the four identities
\begin{align}
& M_{11} M_{12} = M_{12} M_{22}, \hspace{1.55cm}  M_{22} M_{21} = M_{21} M_{11}, \label{eq1} \\
& M_{11} M_{11} - M_{12} M_{21} = M_{11} , \quad M_{22} M_{22} - M_{21} M_{12} = M_{22} .
\label{eq2}
\end{align}
From \eqref{def2M11} and \eqref{def2M12}  it follows  that
\begin{align*}
& M_{11} M_{12} = \begin{bmatrix} T_{+,\a} & T_{+,\ell \b } \end{bmatrix}
\begin{bmatrix} \Delta_{a_0^{-1} } & 0 \\ 0 & - \Delta_{ d_0^{-1} } \end{bmatrix}
\begin{bmatrix} T_{+,\a}^\ast \\ T_{+,\ell \b }^\ast \end{bmatrix} \times\\
&  \hspace{4cm} \times
\begin{bmatrix} H_{+,\ell^\ast \a} & H_{+,\b}  \end{bmatrix}
\begin{bmatrix} - \Delta_{a_0^{-1} } & 0 \\ 0 &  \Delta_{d_0^{-1} } \end{bmatrix}
\begin{bmatrix}  T_{-,\ell^\ast \c}^\ast \\ T_{-,\d }^\ast  \end{bmatrix}.
\end{align*}
Furthermore, from \eqref{formulaM12} and \eqref{def2M22}   it follows  that
\begin{align*}
& M_{12} M_{22} = \begin{bmatrix} T_{+,\a} & T_{+,\ell \b} \end{bmatrix}
\begin{bmatrix} \Delta_{a_0^{-1} } & 0 \\ 0 & - \Delta_{d_0^{-1}} \end{bmatrix}
\begin{bmatrix} H_{-,\c }^\ast \\ H_{-,\ell \d }^\ast \end{bmatrix} \times \\
&\hspace{4cm} \times
\begin{bmatrix} T_{-,\ell^\ast \c} & T_{-,\d}  \end{bmatrix}
\begin{bmatrix} - \Delta_{a_0^{-1} } & 0 \\ 0 & \Delta_{d_0^{-1} } \end{bmatrix}
\begin{bmatrix} T_{-,\ell^\ast \c}^\ast  \\ T_{-,\d}^\ast  \end{bmatrix}.
\end{align*}
But then \eqref{matrixTHHT} shows that $M_{11} M_{12} = M_{12} M_{22}$. In a  similar way, using \eqref{matrixTHHT2} one proves that $ M_{22} M_{21} = M_{21} M_{11} $.

Next, using \eqref{def2M11} and \eqref{formulaM11},  observe that
\begin{align*}
& M_{11} ( M_{11} - I_{ \ell_+^2(\BC^p)}  ) =  -
\begin{bmatrix} T_{+,\a} & T_{+,\ell \b} \end{bmatrix}
\begin{bmatrix} \Delta_{ a_0^{-1} } & 0 \\ 0 & - \Delta_{ d_0^{-1} } \end{bmatrix}
\begin{bmatrix} T_{+,\a }^\ast \\ T_{+,\ell \b}^\ast \end{bmatrix} \times \\
&\hspace{4cm}\times
\begin{bmatrix} H_{+,\ell^\ast \a} & H_{+,\b}  \end{bmatrix}
\begin{bmatrix} - \Delta_{ a_0^{-1} } & 0 \\ 0 &  \Delta_{ d_0^{-1} } \end{bmatrix}
\begin{bmatrix}  H_{+,\ell^\ast \a}^\ast \\ H_{+,\b }^\ast \end{bmatrix}.
\end{align*}
Furthermore, using \eqref{formulaM12} and  \eqref{formulaM21}, we see that
\begin{align*}
& M_{12} M_{21} =
\begin{bmatrix} T_{+,\a} & T_{+,\ell \b} \end{bmatrix}
\begin{bmatrix} \Delta_{ a_0^{-1} } & 0 \\ 0 & - \Delta_{ d_0^{-1} } \end{bmatrix}
\begin{bmatrix}  H_{-,\c }^\ast \\ H_{- ,\ell \d }^\ast \end{bmatrix} \times   \\
&\hspace{4cm} \times
\begin{bmatrix} T_{-,\ell^\ast \c} & T_{-,\d} \end{bmatrix}
\begin{bmatrix} - \Delta_{ a_0^{-1} } & 0 \\ 0 & \Delta_{ d_0^{-1} } \end{bmatrix}
\begin{bmatrix}   H_{+,\ell^\ast \a}^\ast \\ H_{+,\b }^\ast \end{bmatrix}.
\end{align*}
But then \eqref{matrixTHHT} implies that  $M_{11} M_{11} - M_{12} M_{21} = M_{11}$. Similarly, using \eqref{matrixTHHT2} one proves that $ M_{22} M_{22} - M_{21} M_{12} = M_{22}$.

Given  \eqref{invform1} the equalities in \eqref{invform2} are immediate from \eqref{equality1}.
Finally, by multiplying \eqref{rhs3} from the left and the right by
\[
\begin{bmatrix} T_{+,\a} \Delta_{a_0^{-1}} & T_{+,\ell \b } \Delta_{d_0^{-1}}  \end{bmatrix} \ands
\begin{bmatrix} T_{-,\d}^\ast \\ T_{-,\ell^\ast \c}^\ast \end{bmatrix},
\]
respectively, one obtains the equality \eqref{intertwiM}.
\epr

\setcounter{equation}{0}
\section{Direct proof of Theorem \ref{thm:inversionDS2}}\label{sec:prThm3.1}
Let $\a$, $\b$, $\c$, $\d$ be the functions given by \eqref{Fabcd2}, and let $a, b, c, d$ be the associate linear maps given by formulas \eqref{alm}--\eqref{dlm} with both matrices $a_0$ and $d_0$ invertible.
Throughout this section  $ g $ is a function in  $ \sW_+^{p\ts q}$.

\medskip
\noindent
\textbf{Proof of Theorem \ref{thm:inversionDS2}.}
We split the proof  into two parts.

\smallskip
\noindent\textsc{Part 1.} {First we assume that $g$ is a solution of  the twofold EG inverse
problem associated with $ \{ \a,\b , \c ,\d \} $, that is,
we assume that the inclusions in \eqref{incluD12} and \eqref{incluD34} are satisfied.
Our aim is to prove the identity
\begin{equation}\label{inverseOmega}
\begin{bmatrix} M_{11} & M_{12} \\ \noalign{\smallskip} M_{21} & M_{22} \end{bmatrix}
\begin{bmatrix} I_{\ell^2_+ ( \BC^p ) } & H_{+,g} \\ \noalign{\smallskip}
H_{-,g^\ast} & I_{\ell^2_- ( \BC^q ) } \end{bmatrix} =
\begin{bmatrix} I_{\ell^2_+ ( \BC^p ) } & 0 \\ \noalign{\smallskip}
0 & I_{\ell^2_- ( \BC^q ) }  \end{bmatrix}.
\end{equation}
According to Proposition \ref{prop:basicids7} our assumptions imply
that  the operator identities in \eqref{basicids1a}  and \eqref{basicids2a}  hold true.
Furthermore, since $g$ is a solution of  the twofold EG inverse problem associated with
$ \{ \a,\b , \c ,\d \} $, the identities  in \eqref {condD1} (and \eqref {condD1a}) are satisfied, and hence we may use Lemma \ref{lem:altMij}.
We proceed in five steps.

\smallskip
\noindent\textsc{Step 1.1.} According to \eqref{formulaM11}  and \eqref{def2M12} we have
\begin{align*}
& M_{11} - I_{\ell_+^2(\BC^p) } + M_{12} H_{-,g^\ast}=
\\ & =
H_{+,\b} \Delta_{d_0^{-1} } \left( T_{-,\d}^\ast H_{-,g^\ast} + H_{+,\b}^\ast \right)
- H_{+,\ell^\ast \a } \Delta_{a_0^{-1}}
\left( T_{-,\ell^\ast \c}^* H_{-,g^\ast} + H_{+,\ell^\ast \a}^\ast \right).
\end{align*}
Taking the   adjoints of the second identity in \eqref{basicids2a} and of the first identity in \eqref{basicids1a} (using  identities from \eqref{adjointsTH} when necessary), we see that
\[
T_{-,\d}^\ast H_{-,g^\ast} + H_{+,\b}^\ast=0 \ands  T_{-,\ell^\ast \c}^* H_{-,g^\ast} + H_{+,\ell^\ast \a}^\ast =0.
\]
It follows that  $ M_{11} + M_{12} H_{-,g^\ast}  = I_{\ell_+^2(\BC^p)}$.

\smallskip
\noindent\textsc{Step  1.2.}
According to \eqref{def2M11} and \eqref{formulaM12} one has that
\begin{align*}
& M_{11} H_{+,g} + M_{12}= \\
& =
T_{+,\a} \Delta_{ a_0^{-1} }  \left( T_{+,\a}^\ast H_{+,g} + H_{-,\c }^\ast \right)
- T_{+,\ell \b} \Delta_{d_0^{-1} } \left( T_{+,\ell \b}^\ast H_{+,g} + H_{-,\ell \d}^\ast \right).
\end{align*}
Using the second identity \eqref{basicids1a} and the first in \eqref{basicids2a} (together with the identities in \eqref{adjointsTH}) we see that
\begin{equation}\label{2identities}
T_{+,\a}^\ast H_{+,g} + H_{-,\c }^\ast =0\ands T_{+,\ell \b}^\ast H_{+,g} + H_{-,\ell \d}^\ast =0.
\end{equation}
It follows that $ M_{11} H_{+,g} + M_{12} =0 $.

\smallskip
\noindent\textsc{Step  1.3.}
According to \eqref{formulaM22} and \eqref{def2M21} we have
\begin{align*}
& M_{21} H_{+,g} + M_{22} - I_{\ell^2_- ( \BC^q ) } =
\\ &=
H_{-,\c} \Delta_{ a_0^{-1} } \left( H_{-,\c}^\ast + T_{+,\a}^\ast H_{+,g}  \right)
- H_{-,\ell \d} \Delta_{d_0^{-1} } \left( H_{-,\ell \d}^\ast + T_{+,\ell \b }^\ast H_{+,g} \right).
\end{align*}
But then, using the two identities in  \eqref{2identities},  we have
 \[
H_{-,\c}^\ast + T_{+,\a}^\ast H_{+,g}=0\ands H_{-,\ell \d}^\ast + T_{+,\ell \b }^\ast H_{+,g}=0.
\]
It follows that  $  M_{21} H_{+,g} + M_{22} = I_{\ell^2_- ( \BC^q ) } $.

\smallskip
\noindent\textsc{Step  1.4.}
According to  \eqref{formulaM21} and \eqref{def2M22} we have
\begin{align*}
& M_{21}  + M_{22} H_{-,g^\ast } =\\ & =
T_{-,\d} \Delta_{d_0^{-1} } \left( T_{-\d}^\ast H_{-,g^\ast } + H_{+,\b}^\ast \right)
- T_{-,\ell^\ast \c } \Delta_{a_0^{-1} }
\left( T_{-, \ell^\ast \c }^\ast H_{-,g^\ast} + H_{+,\ell^\ast \a}^\ast \right).
\end{align*}
Taking adjoints of the second identity in   \eqref{basicids2a} and of  the first identity in  \eqref{basicids1a} (using  identities from \eqref{adjointsTH} when necessary) we see that
\[
T_{-\d}^\ast H_{-,g^\ast } + H_{+,\b}^\ast =0 \ands T_{-, \ell^\ast \c }^\ast H_{-,g^\ast} + H_{+,\ell^\ast \a}^\ast=0.
\]
It follows that  $ M_{21}  + M_{22} H_{-,g^\ast } = 0 $.

\smallskip
\noindent\textsc{Step  1.5.}
Putting together the results in the preceding four steps  we have
proved \eqref{inverseOmega}.
Since both factors in the left hand side of \eqref{inverseOmega} are selfadjoint,
we have proved \eqref{formO-1}.
Furthermore, since the operator matrix $M$ in \eqref{formO-1} only depends on the given data set $ \{ \a,\b,\c,\d \} $},  so does $\Omega=M^{-1}$ in \eqref{defM1b}. Hence $H_{+,g}$ is uniquely determined by $ \{ \a,\b,\c,\d \} $, and, consequently, it follows that the solution $g$ is unique.

\smallskip
\noindent\textsc{Part 2.}
Conversely, let $g\in \sW_+^{p\ts q}$, and assume that the operator $ \Omega $ given by \eqref{defM1b} is invertible.
Then the equations in the right hand sides of \eqref{incluD12a} and \eqref{incluD34a}
have a unique solution.
Given the solution $ \{a,b,c,d  \} $ define  $ \{\a,\b,\c,\d \} $ by \eqref{defFabcd}. The equivalences in \eqref{incluD12a} and \eqref{incluD34a} imply  that $ \{\a,\b,\c,\d \} $ satisfies \eqref{incluD12} and \eqref{incluD34}.
Thus  indeed $ g \in \sW_+^{ p \times q} $ is a solution of the twofold EG inverse problem
associated with $ \{\a,\b,\c,\d \} $.
\epr


\setcounter{equation}{0}
\section{Proof of Theorem \ref{thm:mainthmDS1}}\label{sec:prThm4.1}
Throughout this section  $\a$, $\b$, $\c$, $\d$ are  the functions  given by \eqref{Fabcd2},  and  $a, b, c, d$ are the associate linear maps given by formulas \eqref{alm}--\eqref{dlm}. We assume that both matrices $a_0$ and $d_0$ are invertible.

\medskip
\noindent\textbf{Proof of Theorem \ref{thm:mainthmDS1}.}
We split the proof into two parts.  The first part concerns the necessity of conditions (D1) and (D2) in  Theorem \ref{thm:mainthmDS1}.

\smallskip\noindent
\textsc{Part 1.}
Suppose that the  twofold  EG inverse problem has a solution. Then, as we mentioned in the introduction (see \eqref{condD1}),  condition (D1) follows from \cite[Theorem 1.2]{KvSch13}.  Theorem \ref{thm:inversionDS2} above states that the operator $ M $ given by
\[
M = \begin{bmatrix} M_{11} & M_{12} \\ M_{21} & M_{22}  \end{bmatrix}
\]
is invertible. It follows from Lemma \ref{lem:inverseR} that the
operators $M_{11} $ and $ M_{22} $ are invertible, hence one-to-one.  In particular, condition (D2) is satisfied.

\smallskip\noindent
\textsc{Part 2.}
In this part we assume that  conditions (D1) and (D2) are satisfied. We show that  $M_{11}$ and $M_{22}$ are invertible, and we prove the reverse implications. This will be done in five steps.

\smallskip\noindent
\textsc{Step 2.1.}
We show that $M_{11}$ and $M_{22}$ are invertible as operators on
$\ell^2$-spaces as well as operators on $\ell^1$-spaces.
Since $\a$, $\b$, $\c$, $\d$ are  Wiener class functions, the corresponding Hankel operators are compact, and hence  \eqref{formulaM11} and  \eqref{formulaM22} imply that both  $M_{11}$ and $M_{22}$
are of the form $ I -K $ with $ I $  an identity operator  and $ K$ a compact operator.
Thus both $ M_{11} $ and $ M_{22} $   are Fredholm operators of index zero, and hence condition (D2) tells us that  these operators are invertible.

Recall that $ \ell_+^1(\BC^p) $  is  contained in    $ \ell_+^2(\BC^p) $ and $ \ell_-^1 (\BC^q) $   in    $ \ell_-^2 (\BC^q) $.
The fact that $\a$, $\b$, $\c$, $\d$ are  Wiener class functions implies that $ M_{11} $ maps
$ \ell_+^1(\BC^p) $ into itself  and $ M_{22}$ maps $ \ell_-^1 (\BC^q) $ into itself.
Thus, by condition (D2), the induced operators
\begin{equation}\label{l1ops}
\tilde{M}_{11} : \ell_+^1 (\BC^p) \to \ell_+^1 (\BC^p) \ands \tilde{M}_{22} : \ell_-^1 (\BC^q) \to \ell_-^1(\BC^q)
\end{equation}
are also one-to-one.
Moreover, again using that   $\a$, $\b$, $\c$, $\d$ are  Wiener class functions, the induced operators  $\tilde{M}_{11}$ and $\tilde{M}_{22}$ are of the form identity minus a compact operator. Hence  these operators  are also invertible. In particular, $M_{11}^{-1}$ maps $ \ell_+^1(\BC^p) $ into itself  and $ M_{22}^{-1}$ maps $ \ell_-^1 (\BC^q) $ into itself.

Finally, with $ M_{11} $ and $ M_{22} $ invertible we obtain from Lemma \ref{lem:inverseR} that $M$ is invertible and
\begin{equation}\label{invM7a}
M^{-1} =
\begin{bmatrix} I_{\ell_+^2(\BC^p) } & -M_{12} M_{22}^{-1} \\ -M_{21} M_{11}^{-1} &
I_{\ell_-^2(\BC^q) } \end{bmatrix}
=
\begin{bmatrix} I_{\ell_+^2(\BC^p) } & -M_{11}^{-1} M_{12} \\ -M_{22}^{-1} M_{21} &
I_{\ell_-^2(\BC^q)}  \end{bmatrix}.
\end{equation}

\smallskip\noindent
\textsc{Step 2.2.}
The next step is to show that $ M_{11}^{-1} M_{12} $ and $ M_{22}^{-1} M_{21} $
are Hankel operators. From \eqref{intertwiM} we know that
$ M_{12} S_{-,q} M_{22} = M_{11} S_{+,p}^\ast M_{12} $.
Therefore
\[
M_{11}^{-1} M_{12} S_{-,q} = S_{+,p}^\ast M_{12} M_{22}^{-1} = S_{+,p}^\ast M_{11}^{-1} M_{12},
\]
with the last identity following from the equality of the right upper corners in \eqref{invM7a}. This intertwining relation  proves that $ M_{11}^{-1} M_{12} $ is a Hankel operator. Since $ M_{11}^\ast = M_{11} $, $ M_{22}^\ast = M_{22} $, and
$ M_{12}^\ast = M_{21} $, we have that
\[
\left( M_{11}^{-1} M_{12} \right)^\ast = M_{21} M_{11}^{-1} = M_{22}^{-1} M_{21}.
\]
It follows that $  M_{22}^{-1} M_{21} $ is also a Hankel operator.

\smallskip\noindent
\textsc{Step 2.3.}
Let $g$ be defined  by the first identity in \eqref{unique1}.
We shall show that $g\in \sW_+^{p\ts q}$ and that  $H_{+, g}=-M_{11}^{-1} M_{12}$.
To do this,  put $h= - M_{11}^{-1} M_{12}\vp_{-,q}$.
From the second identity in \eqref{unitsM} we know that $M_{12}\vp_{-,q}=b$.
Recall  that $ b$ is a linear map from $ \BC^q $ to $ \ell_+^1 (\BC^p) $.
As we have seen in Step 2.1  the operator $ M_{11}^{-1} $  maps   $ \ell_+^1 (\BC^p) $ into $ \ell_+^1 (\BC^p) $.
Hence
\[
h=-M_{11}^{-1} b : \BC^q \to \ell_+^1 (\BC^p)  \ands  g=\sF h\in \sW_+^{p\ts q}.
\]
Since $M_{11}^{-1} M_{12}$ is a Hankel operator, the identities   $h= - M_{11}^{-1} M_{12}\vp_{-,q}$ and $g=\sF h$ imply that $-M_{11}^{-1} M_{12}=H_{+,g}$.

\smallskip\noindent
\textsc{Step 2.4.}
In this part we prove the second identity in \eqref{unique1}.
Put $\tilde{h}=-M_{22}^{-1}M_{21}\vp_{+,p}$.
From the third  identity in \eqref{unitsM} we know that $M_{21}\vp_{+,p}=c$.
Repeating the argument of the previous step, with $c$ in place of $b$ and $M_{22}$ in place of $M_{11}$ we see that
\[
\tilde{h}=-M_{22}^{-1}c:\BC^p\to \ell_-^1(\BC^q) \ands  \tilde{g}:=\sF \tilde{h}\in \sW_-^{q\ts p}.
\]
But then, since $M_{22}^{-1} M_{21}$ is a Hankel operator, we conclude that
$ M_{22}^{-1} M_{21} = - H_{-, \tilde{g}}$, where $\tilde{g}:=\sF \tilde{h}=- \sF( M_{22}^{-1}c)$.
Summarizing, using the results of the two previous steps, we have
\[
H_{-, \tilde{g}}=-M_{22}^{-1} M_{21}=-\left( M_{11}^{-1} M_{12} \right)^\ast =H_{+,g}^\ast=H_{-,g^*}.
\]
Thus $g^*=\tilde{g}=\sF(-M_{22}^{-1}c)$, and the second identity in \eqref{unique1} is proved.

\smallskip\noindent
\textsc{Step 2.5.} \
Finally, we show that $g$ is a solution to the twofold EG inverse problem associated with the data $\{\a,\b,\c,\d \} $.
By Lemma \ref{lem:unitsM} we have
\[
M\mat{c}{\vp_{+,p}\\ 0}=\mat{c}{a\\c}\quad \mbox{and}\quad
M\mat{c}{0\\\vp_{-,q}}=\mat{c}{b\\d}.
\]
Since $M_{11}$ and $M_{22}$ are invertible, the results of the previous steps and the second part of formula \eqref{invform2} show that
\[
\om:= M^{-1}=\begin{bmatrix} I_{\ell_+^2(\BC^p) } & -M_{11}^{-1} M_{12} \\ -M_{22}^{-1} M_{21} &
I_{\ell_-^2(\BC^q)}  \end{bmatrix}= \begin{bmatrix} I_{\ell_+^2(\BC^p) } &H_{+,g}\\ H_{-,g^*}&
I_{\ell_-^2(\BC^q)}  \end{bmatrix}.
\]
We obtain  that
\[
\Omega  \mat{c}{a\\c}=\Omega  M \mat{c}{\vp_{+,p}\\ 0}=\mat{c}{\vp_{+,p}\\ 0},
\]
and similarly
\[
\Omega \mat{c}{b\\d}=\Omega  M \mat{c}{0\\\vp_{-,q}}=\mat{c}{0\\\vp_{-,q}}.
\]
Note that $H_{+,g}=G$, as in \eqref{defG}, and $H_{-,g^*}=G^*$. Hence, via the implications \eqref{incluD12a} and \eqref{incluD34a} we obtain that $g$ satisfies \eqref{incluD12} and \eqref{incluD34}. This shows $g$ is a solution  of the twofold EG inverse problem associated with the data $ \{ \a,\b,\c,\d \}$.
\epr

\setcounter{equation}{0}
\section{Proof of Theorem \ref{thm:posdef}}\label{sec:Contr}
In order to prove Theorem \ref{thm:posdef} we start with a proposition that
on the one hand covers part of Theorem \ref{thm:posdef}, but on the other hand is more detailed.
To state this proposition we need some additional notation.
Recall that $ \sW^{n \times m} $  decomposes as
\[
\sW^{n \times m } = \sW_+^{n \times m } \dot{+} \sW_{-,0}^{n \times m }, \ands
\sW^{n \times m } = \sW_{+,0}^{n \times m } \dot{+}  \sW_-^{n \times m }.
\]
Now let $ \rho \in \sW^{ n \times m} $.
Then, using the above decompositions, we can write $ \rho $
in a unique way as
\[
\rho = \rho_+ + \rho_{-,0}, \ands \rho = \rho_{+,0} + \rho_+,
\]
where $ \rho_+ \in \sW_+^{n \times m }$, $ \rho_{-,0}  \in \sW_{-,0}^{n \times m } $,
$ \rho_{+,0} \in \sW_{+,0}^{n \times m } $ and $ \rho_- \in \sW_-^{n \times m } $.
These direct sum decompositions will play a role in the next proposition.

\begin{prop}\label{prop:posdef}
Let $\{a, b, c, d\}$ be the  data given by formulas \eqref{alm}-- \eqref{dlm},
and let $\a, \b, \c$, and $  \d$ be the functions defined by \eqref{defFabcd}.
\begin{itemize}
\item[\textup{(i)}]
If $ a_0$ is positive definite, $ \a^\ast \a - \c^\ast \c =a_0 $
and $ \det \a(\lambda) $ has no zero in $ \lambda \leq 1 $, then
$ g_1 = - \left( \a^{-\ast} \c^\ast \right)_+ $ is the unique element in $ \sW_+^{p \times q} $
that satisfies the inclusions in \eqref{incluD12}.
\item[\textup{(ii)}]
If $ d_0$ is positive definite, $ \d^\ast \d - \b^\ast \b = d_0 $
and $ \det \d(\lambda) $ has no zero in $ \lambda \geq 1 $, then
$ g_2 = - \left( \b \d^{-1} \right)_+ $ is the unique element in $ \sW_+^{p \times q} $
that satisfies the inclusions in \eqref{incluD34}.
\end{itemize}
Moreover, if in addition the third condition in \eqref{condD1} is also satisfied,
that is $ \a^\ast \b =  \c^\ast \d $, then $ g_1 $ in item \textup{(i)} and $ g_2 $ in item \textup{(ii)} are equal.
\end{prop}

\bpr
It follows from Theorem 3.1 in \cite{KvSch14} that the Fourier coefficients of the unique
function $ g_1 $ that satisfies the inclusions in \eqref{incluD12}
are given by $ - \left( T_{\a^\ast} \right)^{-1} c^\ast $. Therefore
$ g_1 = - \left( \a^{-\ast} \c^\ast \right)_+ $. This proves the first item.

The second item we derive from the first as follows.
Put
\[
\tilde\a(\lambda) = \d(\lambda^{-1} ) , \quad
\tilde\c(\lambda) = \b(\lambda^{-1} ) , \quad
\tilde{p} = q , \quad \tilde{q} = p  \ands \tilde{a}_0 = d_0.
\]
Then item (i) gives that $ \tilde{g}_2 = -  \left( \tilde\a^{-\ast} \tilde\c^\ast  \right)_+ $
is the unique element in $\sW_+^{q \times p } $ such that
\[
\tilde{\a} + \tilde{g}_2 \tilde\c - e_q \in \sW_{-,0}^{q \times q} , \quad
\tilde{g}_2^\ast \tilde\a + \tilde\c \in \sW_{+,0}^{q \times p} .
\]
Put $  g_2 (\lambda ) = \tilde{g}_2( \lambda^{-1} )^\ast $.
Then $ g_2 \in \sW_+^{p \times q} $ and
\[
\d + g_2^\ast \b - e_q \in \sW_{+,0}^{q \times q} , \quad
g_2 \d + \b \in \sW_{-,0}^{p \times q} .
\]
To finish the proof of item (ii) notice that
\begin{align*}
g_2(\lambda)  & = \tilde{g}_2( \lambda^{-1} )^\ast =
- \left[ \left( \tilde\a \bigl( \textstyle{\frac{1}{\lambda}} \bigr)^{-\ast}
\tilde\c \bigl( \textstyle{\frac{1}{\lambda}} \bigr)^{\ast} \right)_- \right]^\ast =
- \left[ \left( \d(\lambda)^{-\ast} \b(\lambda)^\ast \right)_- \right]^\ast \\ & =
-( \b \d^{-1} )_+ (\lambda)  .
\end{align*}

Finally, notice that the conditions in (i) and (ii) imply that
$ \a^{-\ast} $ and $ \d^{-1} $ are well defined and hence in that case $ \a^\ast \b = \c^\ast \d $
implies  $ \a^{-\ast} \c^\ast = \b \d^{-1} $. Therefore
$ \a^\ast \b = \c^\ast \d $ implies $ g_1 = g_2 $.
\epr
\medskip

In the following proof we refer to results in Section \ref{sec:a0d0} that are basically
taken from Sections 3 and 4 in \cite{EGL95}.

\medskip\noindent
\textbf{Proof of Theorem \ref{thm:posdef}.}
Assume that $ g $ is a solution of the EG inverse problem with $ G := H_{+,g} $ strictly contractive.
According to Lemma \ref{lem:a0d0pos12} the matrices $ a_0 $ and $ d_0$ are positive definite.
We know from Theorem 1.2 in \cite{KvSch13}
that condition (ii) is satisfied.
It follows from Proposition \ref{prop:ompos} that $ \det \a $ has no zeros in the unit disk
and $ \det \d(\l) $ has no zeros with $ | \lambda | \geq 1 $.

Conversely, assume that the conditions (i), (ii), and (iii) are satisfied.
If the EG inverse problem has a solution, then Proposition \ref{prop:ompos} tells us that
the operator $ \om_1 $ defined in \eqref{defOmega}
is positive definite, and hence by Lemma \ref{lem:a0d0pos12} we conclude that
$ G $ is strictly contractive.
So it remains to show that there exists a solution.
According to Proposition \ref{prop:posdef} the function $ g_1 = - (\a^{-\ast} \c^\ast)_+ $
is the unique function in $ \sW_+^{p \times q} $ that satisfies the inclusions in \eqref{incluD12}.
Also $ g_2 =  - ( \b \d^{-1})_+ $ is the unique function in $ \sW_+^{p \times q} $
that satisfies the inclusions in  \eqref{incluD34}.
The third statement in Proposition \ref{prop:posdef} gives that $ g_1 = g_2 $.
Hence $ g = g_1 = g_2 $ is the unique solution of the twofold EG inverse problem.
\epr


\setcounter{equation}{0}
\section{The polynomial case.}\label{sec:Polynomial}
In this section we treat the case where the functions $\a$ and $\b$ are polynomials in $ \lambda$,
and $\c$ and $\d$ are polynomials in $ \lambda^{-1}$. We will prove the following theorem.

\begin{thm}\label{thm:polynomial}
Let $\a$ and $\b$ be matrix polynomials in $ \lambda $,
let $\c$ and $\d$ be matrix polynomials in $ \lambda^{-1}$,
and let $m$ be an upper bound of the degrees of $\a$, $\b$, $\c$ and $\d$.
Assume $a_0$ and $d_0$ are invertible.
Then there exists a solution $g$ to the twofold EG-inverse problem
associated with $\a$, $\b$, $\c$, $\d$ if any only if the identities in \eqref{condD1}
are satisfied.
Moreover, this solution $g$ is unique, it is a matrix polynomial
with $\deg g \leq m $ and its coefficients $ g_0, \ldots , g_m  $ are given by
\[
-\mat{ccc}{b_m &\cdots& b_0\\ &\ddots&\vdots\\ &&b_m}\mat{c}{ e_m \\ \vdots\\ e_0 }
=\mat{c}{g_0\\ \vdots\\ g_m}
=-\mat{ccc}{a_0^* &\cdots& a_m^*\\ &\ddots&\vdots\\ &&a_0^*}^{-1}\mat{c}{c_0^*\\ \vdots\\ c_{-m}^*}.
\]
Here
\begin{equation}\label{defe0-m}
\begin{bmatrix} e_m \\ \vdots \\ e_0 \end{bmatrix} :=
\begin{bmatrix} d_0  & \cdots & d_{-m} \\ & \ddots & \vdots  \\ & & d_0   \end{bmatrix}^{-1}
\begin{bmatrix} 0 \\ \vdots \\ 0 \\ I_m \end{bmatrix}.
\end{equation}
\end{thm}

This theorem is closely related to Theorem 2.4 in \cite{EGL96LAA}.
To see this note that under the conditions mentioned in Theorem \ref{thm:polynomial} the equations
\[
\begin{bmatrix} I & G \\ G^\ast & I \end{bmatrix} \begin{bmatrix}  a \\ c  \end{bmatrix} =
 \begin{bmatrix} e_{+,p} \\ 0 \end{bmatrix}, \quad
\begin{bmatrix} I & G \\ G^\ast & I \end{bmatrix}  \begin{bmatrix} b \\ d \end{bmatrix} = \begin{bmatrix} 0 \\  e_{-,q} \end{bmatrix}
\]
are equivalent to
\[
\begin{bmatrix} I & G_m \\ G_m^\ast & I \end{bmatrix} \begin{bmatrix}  a(m) \\ c(m) \end{bmatrix} =
 \begin{bmatrix} e_{+,p}(m) \\ 0 \end{bmatrix}, \quad
\begin{bmatrix} I & G_m \\ G_m^\ast & I \end{bmatrix}  \begin{bmatrix} b(m) \\ d(m) \end{bmatrix} = \begin{bmatrix} 0 \\ e_{-,q}(m) \end{bmatrix},
\]
where $m$ is the upper bound for the degrees given above.
In these equations $G_m$ denotes the right upper $ (m+1) \times (m+1) $ submatrix of $G$,
$ a(m) = \begin{bmatrix} a_0 & \cdots & a_m \end{bmatrix}^\top $,
$ b(m) = \begin{bmatrix} b_0 & \cdots & b_m \end{bmatrix}^\top $,
$ c(m) = \begin{bmatrix} c_{-m} & \cdots & c_0 \end{bmatrix}^\top $ and
$ d(m) = \begin{bmatrix} d_{-m} & \cdots & d_0 \end{bmatrix}^\top $.
The symbol $ e_{+,p}(m) $ denotes the first column of the $ (m+1) \times (m+1) $ identity $ p \times p $ block matrix
and $ e_{-,q}(m) $ denotes the last column of the $ (m+1) \times (m+1) $ identity $ q \times q $ block matrix.
The latter set of equations is solved in \cite[Theorem 2.4]{EGL96LAA} for the case
when $ p=q $.

Here we present a proof of Theorem \ref{thm:polynomial} using Theorem 3.3 in \cite{KvSch14}.
Note that the uniqueness of the solution is covered by Theorem \ref{thm:inversionDS2} above.

We derive Theorem \ref{thm:polynomial} as a corollary of the following proposition.

\begin{prop}\label{prop:polynomial}
Let $\a$ and $\b$ be matrix polynomials in $\lambda $,
let $\c$ and $\d$ be matrix polynomials in $ \lambda^{-1}$,
and let $m$ be an upper bound of the degrees of $\a$, $\b$, $\c$ and $\d$.
Assume $a_0$ and $d_0$ are invertible.
\begin{itemize}
\item[\textup{(i)}]
If $ \a^\ast \a - \c^\ast \c =a_0 $, then the onefold EG inverse problem associated
with $ \a $ and $ \c $ has a unique polynomial solution $ g $, which has degree at most $ m $.
Moreover the Fourier coefficients $ g_0, \ldots , g_m $ of $g$ are given by
\[
\mat{c}{g_0\\ \vdots\\ g_m}
=-\mat{ccc}{a_0^* &\cdots& a_m^*\\ &\ddots&\vdots\\ &&a_0^*}^{-1}\mat{c}{c_0^*\\ \vdots\\ c_{-m}^*}.
\]
\item[\textup{(ii)}]
If  $ \d^\ast \d - \b^\ast \b = d_0 $, then there exists a unique polynomial
$ \varphi (\lambda) = \sum_{j=0}^m \varphi_{j} \lambda^{-j} $ in
$ \sW_-^{q \times p} $ such that
\[
\d + \varphi  \b  - e_q \in \sW_{+,0}^{q \times q} \ands
\varphi^\ast \d + \b \in \sW_{-,0}^{p \times q}.
\]
Moreover the Fourier coefficients $ \varphi_0 , \ldots, \varphi_{-m} $ are given by
\begin{align}
& \begin{bmatrix} \varphi_0 & \cdots & \varphi_{m} \end{bmatrix}  = \nonumber \\
& \qquad \ - \begin{bmatrix} 0 & \cdots & 0  & I_q \end{bmatrix}
\begin{bmatrix} d_0^* & & \\ \vdots & \ddots &  \\ d_{-m}^* & \cdots & d_0^*  \end{bmatrix}^{-1}
\begin{bmatrix} b_m^* && \\ \vdots & \ddots & \\ b_0^* & \cdots & b_m^* \end{bmatrix}.
\label{EGdb}
\end{align}
\end{itemize}
Moreover, if all three conditions in \eqref{condD1} are satisfied, then for $ g $ and $ \varphi $
as in items \textup{(i)} and \textup{(ii)} one has  $ \varphi^\ast = g $.
\end{prop}

The first statement of this proposition can be found in \cite[Theorem 2.1]{EGL96LAA}
for the case when $ p= q$.

Before we start the proof we introduce some notation.
We denote the compressions to the first $ m+1 $ components
of $ T_{+,\a} $ by $ T_{+,\a,m} $,  of $ H_{+,\b } $ by $ H_{+,\b,m} $,
of $ H_{- ,\c } $ by $ H_{-,\c,m} $, and of $ T_{-,\d} $ by $ T_{- ,\d,m} $.
Thus  $ y = T_{+,\a,m} x $ and $ y = H_{+,\b,m} x $ if and only if
\[
\begin{bmatrix} y_0 \\ \vdots \\ y_m  \end{bmatrix} =
\begin{bmatrix} a_0  & & \\ \vdots & \ddots & \\ a_m & \cdots & a_0  \end{bmatrix}
\begin{bmatrix} x_0 \\ \vdots \\ x_m  \end{bmatrix},\quad \begin{bmatrix} y_0 \\ \vdots \\ y_m  \end{bmatrix} =
\begin{bmatrix} b_m  & \cdots & b_0 \\ & \ddots & \vdots  \\ & & b_m   \end{bmatrix}
\begin{bmatrix} x_{-m} \\ \vdots \\ x_0  \end{bmatrix},
\]
respectively. Similarly, $ y = H_{-,\c,m} x $  and  $ y = T_{-,\d,m} x $ if and only if
\[
\begin{bmatrix} y_{-m} \\ \vdots \\ y_0  \end{bmatrix} =
\begin{bmatrix} c_{-m}  & & \\ \vdots & \ddots & \\ c_0 & \cdots & c_{-m} \end{bmatrix}
\begin{bmatrix} x_0 \\ \vdots \\ x_m  \end{bmatrix}, \quad
\begin{bmatrix} y_{-m} \\ \vdots \\ y_0  \end{bmatrix} =
\begin{bmatrix} d_0  & \cdots & d_{-m} \\ & \ddots & \vdots  \\ & & d_0   \end{bmatrix}
\begin{bmatrix} x_{-m} \\ \vdots \\ x_0  \end{bmatrix},
\]
respectively.

\medskip
\textbf{Proof of Proposition \ref{prop:polynomial}.}
Item (i) is a direct consequence of Theorem 3.3 in \cite{KvSch14}.
For later purpose we remark that
\[
\begin{bmatrix} g_0 \\ \vdots \\ g_m  \end{bmatrix} =
- T_{+,\a,m}^{-\ast} H_{-,\c,m}^\ast \begin{bmatrix} 0 \\ \vdots \\0  \\ I_q  \end{bmatrix}.
\]

The next step is to prove item (ii).
Define polynomials $\tilde{\a}( \lambda )=\d(\lambda^{-1})$ and
$\tilde{\c}(\lambda)=\b( \lambda^{-1} )$ in $\lambda$.
Note that $\tilde{\a}(0)=d_0$ is invertible, and that we have that
\[
\tilde{\a}^\ast \tilde{\a} - \tilde{\c}^\ast \c = d_0.
\]
Again applying Theorem 3.3 from \cite{KvSch14} we obtain that there exists
a unique matrix polynomial $ \tilde{\varphi}(\lambda) = \sum_{j=0}^m \varphi_{j} \lambda^j $ so that
\[
\tilde{\a}+ \tilde{\varphi} \tilde{\c} - e_q \in  \sW _{-,0}^{q \times q} \ands
\tilde{\varphi}^\ast \tilde{\a} + \tilde{c} \in  \sW_{+,0}^{p \times q}.
\]
Moreover, this matrix polynomial $ \tilde{\varphi} $ satisfies $\deg \tilde{\varphi} \leq m$
and its coefficients are given by
\begin{equation*} 
\mat{c}{ \varphi_0 \\ \vdots\\ \varphi_m }
= - \mat{ccc}{d_0^* &\cdots& d_{-m}^*\\ &\ddots&\vdots\\ &&d_0^*}^{-1}
\mat{c}{b_0^*\\ \vdots\\ b_{m}^*}  ,
\end{equation*}
or equivalently
\begin{align}
\begin{bmatrix} \varphi_0 & \cdots & \varphi_{m} \end{bmatrix} & =  -
\begin{bmatrix} 0 & \cdots & 0  & I_q \end{bmatrix}
\begin{bmatrix} d_0^* & & \\ \vdots & \ddots &  \\ d_{-m}^* & \cdots & d_0^*  \end{bmatrix}^{-1}
\begin{bmatrix} b_m^* && \\ \vdots & \ddots & \\ b_0^* & \cdots & b_m^* \end{bmatrix}
\nonumber \\
& = - e_{-,q}(m)^* T_{-,\d,m}^{-*} H_{+,\b,m}^* . \label{phicoeff}
\end{align}
Here we use that the inverse of an invertible block triangular Toeplitz matrix
is again a block triangular Toeplitz matrix.
So there exist  $ q \times q$ matrices $ e_0 , \ldots , e_m $ such that
\[
\mat{ccc}{d_0^* &\cdots& d_{-m}^*\\ &\ddots&\vdots\\ &&d_0^*}^{-1} =
\begin{bmatrix} e_0^* & \cdots & e_m^* \\ & \ddots & \vdots \\ && e_0^* \end{bmatrix},
\]
and hence
\[
\begin{bmatrix} d_0^* & & \\ \vdots & \ddots &  \\ d_{-m}^* & \cdots & d_0^*  \end{bmatrix}^{-1} =
\begin{bmatrix} e_0^* & & \\ \vdots & \ddots & \\ e_m^* & \cdots & e_0^* \end{bmatrix}.
\]
Put $ \varphi(\lambda) = \tilde{\varphi} (\lambda^{-1} ) $.
Then
\[
\d + \varphi \b - e_q \in \sW_{+,0}^{q \times q} \ands
\varphi^* \d + \b \in  \sW_{-,0}^{p \times q}.
\]

To finish the proof we need to show that $ \varphi ( \lambda )^\ast = g( \lambda ) $.
In other words we need to prove that
\[
\psi^* = \begin{bmatrix} \varphi_0^* \\ \vdots \\ \varphi_{m}^\ast  \end{bmatrix} =
\begin{bmatrix} g_0 \\ \vdots \\ g_m \end{bmatrix} = h .
\]
We claim that this is the case as a result of the third identity in \eqref{condD1}.
As observed in Lemma \ref{basicids5} we have
\[
T_{+,\a}^* H_{+,\b}=T_{+,\a^*} H_{+,\b}=H_{+,\c^*} T_{-,\d}=H_{-,\c}^* T_{-,\d}.
\]
Since
\[
H_{+,\b} = \begin{bmatrix} 0 & H_{+,\b,m} \\ 0 & 0 \end{bmatrix} \ands
H_{-,\c} = \begin{bmatrix} 0 & 0 \\  H_{-,\c,m}  & 0 \end{bmatrix} ,
\]
we have $ T_{+,\a,m}^* H_{+,\b,m}= H_{-,\c,m}^* T_{-,\d,m} $.
It follows that
\[
\psi^* = - H _{+,\b,m} T_{-,\d,m}^{-\ast} e_{-,q}(m) = - T_{+,\a,m}^{-*} H_{-,\c,m}^* e_{-,q}(m) = h.
\]
\epr


\setcounter{equation}{0}
\appendix
\renewcommand{\theequation}{A.\arabic{equation}}
\section{the role of the matrices $a_0$ and $d_0$ }\label{sec:a0d0}
In the main theorems of this paper it is assumed that the matrices
$ a_0 $ and $ d_0$ are invertible.
On the one hand these conditions can be weakened.
For example, in many cases it suffices to assume that only one of the two is invertible.
On the other hand additional conditions on $ a_0 $ and $ d_0 $ yield additional properties of the solution of the EG inverse problem.
These facts can be found in a somewhat less general form in
Sections 3 and 4 in \cite{EGL95}, where the connections between polynomials and finite Toeplitz matrices is a starting point (cf., Chapter 1 in \cite{EG03}). For the sake of completeness these results are reviewed in this section.
We apply them in, e.g., the proof of Proposition \ref{prop:posdef}.

Let $g$ be any function in $\sW_+^{p\ts q}$, and let  $g(\l) = \sum_{\nu=0}^\iy  \l^{\nu} g_\nu$, $\l \in \BT$.  With $g$ we associate the Hankel operator $G$ defined by
\eqref{defG}.
Put
\begin{equation}\label{defG1}
G_1=\begin{bmatrix}
\cdots&g_3&g_2&g_1\\
\cdots& g_4&g_3&g_2\\
\cdots& g_5&g_4&g_3\\
&\vdots&\vdots&\vdots
\end{bmatrix}: \ell_-^2(\BC^q)\to \ell_+^2(\BC^{p}).
\end{equation}
Note that
\[
G_1 = S_{+,p}^* G = G S_{-,q}  ,
\]
where $ S_{+,p} $ and $ S_{-,q}$ are defined by \eqref{shiftminp} and \eqref{shiftplusq}, respectively.
We have
\begin{align}
G&=\begin{bmatrix}R_1\\G_1\end{bmatrix}, \quad \mbox{with } R_1=\begin{bmatrix}\cdots &g_2&g_1&g_0\end{bmatrix}: \ell_-^2(\BC^q)\to \BC^p, \label{partition1a}\\
G&=\begin{bmatrix}G_1& K_1\end{bmatrix}, \quad \mbox{with } K_1= \begin{bmatrix}g_0\\ g_1\\ g_2\\ \vdots \end{bmatrix}:\BC^q\to \ell_+^2(\BC^{p}).\label{partition1b}
\end{align}
Next define
\begin{equation}\label{defOmega}
\om := \begin{bmatrix}I& G\\ G^*&I\end{bmatrix} \ands \om_1 := \begin{bmatrix}I& G_1\\ G_1^*&I\end{bmatrix}.
\end{equation}

\begin{prop}\label{prop:a0d0inv}
Assume that $ g $ is a solution to the twofold EG inverse problem associated
with the date set $ \{ \a, \b, \c, \d  \} $.
\begin{itemize}
\item[\textup{(i)}]
Then  $ a_0 $ and $ d_0$ are selfadjoint.
If $ a_0 $ or $ d_0 $ is invertible, then $ \Omega $ is invertible.
\item[\textup{(ii)}]
Conversely, if $ \Omega $ is invertible, then $ \Omega_1 $ is invertible if and only if
$ a_0 $ or $ d_0 $ is invertible,
and in that case both $ a_0 $ and $ d_0 $ are invertible.
\end{itemize}
\end{prop}

\bpr
First recall that Theorem 1.2 in \cite{KvSch13} gives that \eqref{condD1} is satisfied
and hence $ a_0 $ and $ d_0 $ are selfadjoint.
Theorem 1.1 in \cite{KvSch12} gives that if $ a_0 $ or $ d_0 $ is invertible, then
$ \Omega $ is invertible.
This proves item (i).

In order to prove item (ii) we need some preliminaries.
Using the partitionings in \eqref{partition1a} and \eqref{partition1b} we see that the operator $\om$ admits the following  $3 \ts 3$ block partitionings :
\begin{equation}\label{partition2}
\om = \begin{bmatrix}
I_p & 0 & | & R_1\\
0 & I & | & G_1\\
-- &-- & | & --\\
R_1^* & G_1^ *& | & I
\end{bmatrix} \ands
\om = \begin{bmatrix}
I & | & G_1 & K_1\\
-- & | & -- & -- \\
G_1^* & | & I & 0 \\
K_1^* & | & 0 & I_q
\end{bmatrix}
\end{equation}
Using the definition of $ \om_1 $ in the right hand side of \eqref{defOmega}  we obtain two alternative $ 2\ts 2 $ block operator matrix representations of  $ \om $:
\begin{equation}\label{partition3}
\om=\begin{bmatrix}
I_p&R\\
R^*&\om_1
\end{bmatrix} \ands
\om=\begin{bmatrix}
\om_1&K\\
K^*&I_q
\end{bmatrix}
\end{equation}
Here
\begin{equation}\label{defKR}
R=\begin{bmatrix} 0&R_1 \end{bmatrix}\ands  K=\begin{bmatrix}K_1\\0 \end{bmatrix}.
\end{equation}
Now use that $g\in \sW_+^{p\ts q}$ is a solution to the twofold EG inverse problem associated with the data set  $\{\a, \b, \c, \d\}$.
Thus
\begin{equation}\label{2fold1}
\begin{bmatrix}I&G\\ G^*&I\end{bmatrix}
\begin{bmatrix} a\\ c\end{bmatrix}=\begin{bmatrix} \vp_{+,p}\\0\end{bmatrix}\ands
\begin{bmatrix}I&G\\ G^*&I\end{bmatrix}
\begin{bmatrix} b\\ d\end{bmatrix}=\begin{bmatrix}0\\ \vp_{-,q}\end{bmatrix}
\end{equation}
Using the partitions in \eqref{partition2}, the identities in \eqref{2fold1} can be rewritten as
\begin{equation}\label{2fold2}
\begin{bmatrix}
I_p&R\\
R^*&\om_1
\end{bmatrix} \begin{bmatrix} a_0\\ X \end{bmatrix}=\begin{bmatrix} I_p\\ 0 \end{bmatrix} \ands \begin{bmatrix}
\om_1&K\\
K^*& I_q
\end{bmatrix}\begin{bmatrix} Y\\ d_0 \end{bmatrix}= \begin{bmatrix} 0\\ I_q \end{bmatrix}
\end{equation}

Assume now that $ \Omega $ is invertible.
Then the two identities  in \eqref{2fold2} tell us that
\begin{equation}\label{Schurpr}
\begin{bmatrix} I_p & 0 \end{bmatrix}\om^{-1}\begin{bmatrix} I_p \\ 0 \end{bmatrix}=a_0 \ands \begin{bmatrix} 0 & I_q  \end{bmatrix}\om^{-1}\begin{bmatrix} 0 \\I_q\end{bmatrix}=d_0.
\end{equation}
Moreover, using the two  identities in \eqref{Schurpr},
a standard Schur complement argument shows that
\begin{align*}
&\mbox{$a_0$ is invertible} \quad \Longleftrightarrow \quad  \mbox{$\om_1$ is invertible}, \\
&\mbox{$d_0$ is invertible} \quad \Longleftrightarrow \quad  \mbox{$\om_1$ is invertible}.
\end {align*}
Hence item (ii) is proved.
\epr

\begin{lem}\label{lem:a0d0pos12}
Assume that $g\in \sW_+^{p\ts q}$ is a solution to the twofold EG inverse problem associated with the data set $\{\a, \b, \c,\d\}$.
If $G$ in \eqref{defG} is strictly contractive, then  both $a_0$ and $d_0$ are positive definite and
$ \Omega_1 $ is strictly positive.
Conversely, if $ \om_1 $ is strictly positive and $a_0$ or $d_0$ is  positive definite, then $G$ is strictly contractive.
\end{lem}

\bpr
The two  identities in \eqref{2fold2}  also yield  the following  identities:
\begin{align}
&\begin{bmatrix} a_0&X^*\\ 0 &I   \end{bmatrix}\om \begin{bmatrix} a_0&0\\ X &I   \end{bmatrix}=\begin{bmatrix} a_0&0\\ 0 &\om_1   \end{bmatrix},\label{2fold3a} \\[.2cm]
&\begin{bmatrix}I&0\\ Y^* &d_0   \end{bmatrix}\om \begin{bmatrix} I&Y\\ 0 &d_0   \end{bmatrix}=\begin{bmatrix} \om_1&0\\ 0 &d_0  \end{bmatrix}.\label{2fold3b}
\end {align}
Indeed,
\[
\begin{bmatrix} a_0&X^*\\ 0 &I   \end{bmatrix}\om \begin{bmatrix} a_0&0\\ X &I   \end{bmatrix}=\begin{bmatrix} a_0&X^*\\ 0 &I   \end{bmatrix}\begin{bmatrix} I_p&R\\ 0& \om_1\end{bmatrix}=\begin{bmatrix} a_0&a_0 R+X^*\om_1\\ 0 &\om_1   \end{bmatrix}.
\]
From the first identity in \eqref{2fold2} we know that $R^* a_0 +\om_1 X=0$.
Since both $a_0$ and $\om_1$ are selfadjoint, it follows that  $ a_0 R+ X^* \om_1=0$,
and \eqref{2fold3a} is proved.  The identity \eqref{2fold3b} is proved in a similar way.

Note that $G$ is strictly contractive if and only if $  \om $ is strictly positive.
The identity \eqref{2fold3a} gives that $\om$ is strictly positive if and only if
$ \om_1 $ is strictly positive and $a_0$ is positive definite.
Similarly, the identity \eqref{2fold3b} gives that $\om$ is strictly positive if and only if
$ \om_1 $ is strictly positive and $d_0$ is  positive definite.
\epr

\medskip\noindent
For $ p = q $ the following proposition is an immediate consequence of
Theorems 4.1 and 4.2  in \cite{EGL95}.
If $ p \not= q$, with a minor modification of the data $\a$, $ \b $, $ \c $, $\d$
the arguments used to prove Theorems 4.1 and 4.2 in \cite{EGL95} also yield the result below.

\begin{prop}\label{prop:ompos} Let $g\in \sW_+^{p\ts q}$ be a solution to the twofold EG inverse problem associated with the data set $\{\a, \b, \c, \d\}$ with $ \a, \b, \c$, and $\d$ the functions defined by \eqref{defFabcd}.
\begin{itemize}
\item[\textup{(i)}]  Assume  $a_0$ is positive definite. Then  $\det \a$ has no zero on $\BT$, and $\om_1$ is strictly positive if and only if $\det \a$ has no zero inside  the unit circle.
\item[\textup{(ii)}]  Assume  $d_0$ is positive definite. Then  $\det \d$ has no zero on $\BT$, and  $\om_1$ is strictly positive if and only if $\det \d$ has no zero outside the unit circle.
\end{itemize}
\end{prop}

\bpr
We only have to consider the case when $ p \not= q$.
Assume $ p > q $.
Let $ \a, \b, \c, \d $ be the functions defined by \eqref{defFabcd}.
Put
\[
\tilde{\a} = \a, \quad \tilde{\b} = \begin{bmatrix} \b & 0  \end{bmatrix}, \  
\tilde{\c} = \begin{bmatrix} \c \\ 0  \end{bmatrix}, \  
\tilde{\d} = \begin{bmatrix} \d & 0 \\ 0 & I_{p-q} \end{bmatrix},  \
\tilde{g} =  \begin{bmatrix} g & 0  \end{bmatrix}.
\]
Here $ 0 $ stands for the a zero matrix which each time is chosen such that the extended matrix
is of size $ p \times p $.
Thus $ \tilde\a \in \sW_+^{p \times p } $, $\tilde\b \in \sW_+^{p \times p }$,
$\tilde\c \in \sW_-^{p \times p }$, $ \tilde\d \in \sW_-^{p \times p }$, and
$ \tilde{g} \in \sW_+^{p \times p } $.
Let $  \{ \tilde{a}, \tilde{b}, \tilde{c}, \tilde{d} \} $ be the dataset
corresponding to $ \tilde\a, \tilde\b, \tilde\c$, $ \tilde\d $ as in \eqref{defFabcd}.
Then it follows by direct verification that $ \tilde{g} $ is the solution of the
twofold EG  inverse problem associated to the dataset
$ \{ \tilde\a, \tilde\b, \tilde\c, \tilde\d \} $.
Let $ \tilde{\Omega}_1 $ be defined as $ \Omega_1 $ with $ g $ replaced by $ \tilde{g} $.
Now assume that $ \tilde{a_0} = a_0 $ is positive definite.
Then it follows from Theorem 4.1  in \cite{EGL95} that
$ \det\tilde\a = \det\a $ has no zero on $ \BT$, and $ \tilde{\Omega}_1 $ is strictly
positive if and only if $ \det\tilde\a$ has no zero inside $ \BT $.
Notice that there exists an invertible transformation $ E $ such that
\[
\tilde{\Omega}_1 = E \begin{bmatrix} \Omega_1 & 0 \\ 0 & I \end{bmatrix} E^{-1}.
\]
In particular we get that $ \tilde{\Omega}_1 $ is positive definite
if and only if $ \Omega_1 $ is.

The case when $ p < q $ and the item (ii) are proved in a similar way.
\epr

\medskip

\paragraph{\bf Acknowledgement}
This work is based on the research supported in part by the National
Research Foundation of South Africa. Any opinion, finding and conclusion or
recommendation expressed in this material is that of the authors and
the NRF does not accept any liability in this regard.

\end{document}